\documentclass[11pt,amsfonts]{article}
\usepackage{graphicx,amssymb,eucal,mathrsfs,color}
\usepackage{latexsym,amsmath,epsfig,epic,eepic}
\usepackage[all]{xy}
\usepackage{wasysym,ntheorem}

\newcommand{\Z}{\mathbb{\bf Z}}
\newcommand{\R}{\mathbb{\bf R}}

\newcommand{\C}{\mathbb{C}}

\newcommand{\PSL}{\mathrm{PSL}_{2}(\mathbb{C})}
\newcommand{\PSLR}{\mathrm{PSL}_{2}(\mathbb{R})}

\newcommand{\HH}{\mathbb{H}^{n}}
\newcommand{\HHI}{\mathbb{H}^{\infty}}
\newcommand{\SO}{{\rm Isom}(\mathbb{H}^{n})}
\newcommand{\SU}{{\rm SU}(n,1)}
\newcommand{\bir}{{\rm Bir}(\mathbb{P}^{2})}
\newcommand{\PGL}{{\rm PGL}(3,\mathbb{C})}
\newcommand{\GL}{{\rm GL}(\mathscr{H})}

\newcommand{\ZZZ}{\mathscr{Z}(\mathbb{P}^{2})}
\newcommand{\ZZ}{{\rm Z}(\mathbb{P}^{2})}

\newtheorem{main}{{\bf  Theorem}}
\newtheorem{prop}{Proposition}
\newtheorem{lemma}[prop]{Lemma}
\newtheorem{rem}{Remark}
\newtheorem{cor}{Corollary}
\newtheorem{exam}{Example}
\newtheorem{defi}{Definition}

\title{K\"{a}hler groups, real hyperbolic spaces and the Cremona group}
\author{\textsc Thomas Delzant \and Pierre Py \thanks{Supported in part by NSF grant DMS 0905911}}
\date{July 2011}

\begin{document}
\maketitle

\begin{abstract} Generalizing a classical theorem of Carlson and Toledo, we prove that {\it any} Zariski dense isometric action of a K\"{a}hler group on the real hyperbolic space of dimension at least $3$ factors through a homomorphism onto a cocompact discrete subgroup of $\PSLR$. We also study actions of K\"{a}hler groups on infinite dimensional real hyperbolic spaces, describe some exotic actions of $\PSLR$ on these spaces, and give an application to the study of the Cremona group. 
\end{abstract}

\tableofcontents
\newpage

\section{Introduction}
\subsection{K\"{a}hler groups and real hyperbolic spaces}

In this text, we study actions of fundamental groups of compact K\"{a}hler manifolds (referred to as {\it K\"{a}hler groups}) on finite or infinite dimensional real hyperbolic spaces. For an introduction to the study of K\"{a}hler groups, the reader can consult \cite{abckt} as well as \cite{burger,corsimp,klingler1,kkm,napram,reznikov} for more recent developments. 

All along the text, we denote by $X$ a compact K\"{a}hler manifold and by $\Gamma$ its fundamental group. 

First, it is a classical result of Carlson and Toledo \cite{ct1} that a cocompact lattice of the group $\SO$ of isometries of the real hyperbolic space $\HH$ of dimension $n\ge 3$ cannot be isomorphic to a K\"{a}hler group (see also \cite{abckt}, section 6.4, for a proof in the nonuniform case, as well as \cite{jostyau,yihu1,yihu2}). Their result is actually more precise. Let us first make the following classical definitions.  

\begin{defi}\label{orbifold} A hyperbolic $2$-orbifold $\Sigma$ is a quotient of the unit disc $\Delta$ in $\C$ by a discrete cocompact subgroup (denoted by $\pi_{1}^{orb}(\Sigma)$) of $\PSLR$. A map from a complex manifold $X$ to $\Sigma$ is holomorphic if it lifts to a holomorphic map from the universal cover of $X$ to the unit disc. 
\end{defi}

Note that we could have given an equivalent definition of a hyperbolic $2$-orbifold in terms of Riemann surfaces with marked points (see \cite{delzant2}); however Definition~\ref{orbifold} will be sufficient for our purpose.  

\begin{defi} A fibration of a compact complex manifold $X$ onto a hyperbolic $2$-orbifold $\Sigma$ is a holomorphic surjective map $f : X \to \Sigma$ with connected fibers. 
\end{defi}

We say that a linear representation $\rho : \Gamma \to {\rm GL}_{{\rm N}}(\C)$ factors through a fibration \mbox{$f : X \to \Sigma$} onto a hyperbolic $2$-orbifold if there exists a linear representation $$\hat{\rho} : \pi_{1}^{orb}(\Sigma)\to {\rm GL}_{{\rm N}}(\C)$$ 
\noindent such that $\rho =\hat{\rho}\circ f_{\ast}$, where $f_{\ast} : \Gamma \to \pi_{1}^{orb}(\Sigma)$ is the homomorphism induced by $f$. In \cite{ct1}, Carlson and Toledo, relying on a result of Sampson \cite{sampson2}, proved that any homomorphism from a K\"{a}hler group onto a cocompact lattice in $\SO$ ($n\ge 3$) factors through a fibration as above. Note that Carlson and Toledo's proof still applies and yields the same result for homomorphisms $\Gamma \to \SO$ with Zariski dense and discrete image. We will prove that the same result remains true for any homomorphism into $\SO$ ($n\ge 3$) with Zariski dense image, discrete or not.

We now want to study isometric actions of $\Gamma$ on infinite dimensional real hyperbolic spaces. For an introduction to these spaces, see \cite{gromov1}, section 6, and \cite{bim}. Let $\mathscr{H}$ be a real Hilbert space with scalar product $\langle \cdot , \cdot \rangle$. We always assume that $\mathscr{H}$ is separable. Let $u$ be a unit vector in $\mathscr{H}$, and define a quadratic form $B$ on $\mathscr{H}$ by:
$$B(v_{1}+t_{1}u,v_{2}+t_{2}u)=t_{1}t_{2}-\langle v_{1},v_{2}\rangle ,$$
\noindent where the $v_{i}$ are orthogonal to $u$ and $t_{i}\in \mathbb{R}$. Copying the construction of finite dimensional hyperbolic spaces, one defines the hyperbolic space $\HHI$ associated to $\mathscr{H}$ by: 
$$\HHI:=\{v\in \mathscr{H}, B(v,v)=1,  B(v,u) >0\}.$$

\noindent The formula 
$${\rm cosh} (d(x,y))=B(x,y)\;\;\;\;\; (x,y\in \HHI),$$

\noindent defines a distance on $\HHI$ which turns it into a complete ${\rm CAT}(-1)$ metric space. We denote by ${\rm Isom}(\HHI)$ the group of isometries of $\HHI$ and by $\partial \HHI$ its boundary. We say that an isometric action of a finitely generated group on $\HHI$ is {\it elementary} if it  fixes a geodesic in $\HHI$ or a point in $\HHI \cup \partial \HHI$; {\it non-elementary} otherwise. On the other hand we say that a group of isometries of $\HHI$ is {\it elliptic} if it fixes a point in $\HHI$.  We can then prove:

\begin{main}\label{infinite} Let $\Gamma$ be a K\"{a}hler group. Let $\rho : \Gamma \to {\rm Isom}(\HHI)$ be a non-elementary action of $\Gamma$ on $\HHI$. Assume that $\rho$ is minimal, i.e. that $\HHI$ contains no nontrivial closed $\rho$-invariant totally geodesic subspace. Then, one of the following two cases happens. 
\begin{enumerate}
\item The representation $\rho$ factors through a fibration onto a hyperbolic $2$-orbifold.
\item  The representation $\rho$ can be written as $\rho =\Psi \circ \theta$, where $\theta$ is a homomorphism  from $\Gamma$ to $\PSLR$ with dense image and $\Psi : \PSLR \to {\rm Isom}(\HHI)$ is a continuous homomorphism. 
\end{enumerate}
\end{main}

Let us make some comments about this result. First, according to \cite{bim}, Proposition 4.3, if $\rho : \Gamma \to {\rm Isom}(\HHI)$ is a non-elementary homomorphism, there exists a unique closed totally geodesic subspace $\HHI_{\rho}$ of $\HHI$ which is $\Gamma$-invariant and minimal with respect to this property. Hence if $\rho$ is nonminimal, one can always apply the previous theorem to the induced action on $\HHI_{\rho}$.

Second, it will be clear from the proof that Theorem~\ref{infinite} remains true when $\HHI$ is replaced by a finite dimensional hyperbolic space $\HH$. But the second possibility in Theorem~\ref{infinite} never occurs for finite dimensional spaces (and if $n\ge 3$).  This is a consequence of the fact that any continuous action of $\PSLR$ on $\HH$  has to preserve a totally geodesic plane (see \cite{gromovpansu} page 122, as well as \cite{karpe,mostow}), hence cannot be minimal. Hence, we obtain the:
\begin{cor}
If $n\ge 3$, any homomorphism $\Gamma \to \SO$ with Zariski dense image factors through a fibration onto a hyperbolic $2$-orbifold.  
\end{cor}

Third, let us recall that in \cite{corsimp}, Corlette and Simpson classified Zariski dense representations of K\"{a}hler groups into $\PSL$ (see \cite{delzant2} for the non-Zariski dense case): any such representation either factors through a fibration onto a hyperbolic $2$-orbifold or factors through a holomorphic map to a compact manifold whose universal cover is a product of discs. Observe also that when $n=3$ the group ${\rm Isom}(\mathbb{H}^{3})$ is isomorphic to $\PSL$ (up to a $2$-sheeted cover), and the corollary above contains only one case, as opposed to the result of \cite{corsimp}.  This is because we ask that the {\it real} Zariski closure of $\rho$ is $\PSL$, instead of merely asking that its {\it complex} Zariski closure is $\PSL$.  Observe also that when the second case of Theorem~\ref{infinite} occurs, we can apply the result of Corlette and Simpson to the representation $\theta$.

The key fact we just used (namely, that any action of $\PSLR$ on a finite dimensional hyperbolic space preserves a totally geodesic plane) is not true for actions on $\HHI$: there exist irreducible continuous linear representations  of $\PSLR$ on a Hilbert space, which are not unitary but which preserve a quadratic form of signature $(\infty,1)$. The associated actions on $\HHI$ are minimal, hence do not have an invariant totally geodesic plane. As a consequence, the second possibility from Theorem~\ref{infinite} actually occurs in some examples. Although this might be well-known to experts in representation theory, we will describe these representations in section~\ref{sl2}. There, we will prove the following theorem.

\begin{main}\label{minim} There exists a one parameter family of representations $$\rho_{t} : \PSLR \to {\rm Isom}(\HHI)$$ where $t\in (0,1)$ with the following properties:
\begin{enumerate}
\item The action $\rho_{t}$ is non-elementary and has no nontrivial closed totally geodesic invariant subspace.
\item There exists a $\rho_{t}$-equivariant harmonic map $f_{t} : \Delta \to \HHI$ whose image is a minimal surface with curvature $\frac{-2}{t(t+1)}$. 
\end{enumerate}
\end{main}

\subsection{Cremona group and Picard-Manin space}

The main example of a group acting on an infinite dimensional hyperbolic space comes from algebraic geometry. Recall that the Cremona group, denoted by $\bir$, is the group of birational transformations of the complex projective plane $\mathbb{P}^{2}$. A birational transformation of $\mathbb{P}^{2}$ is a transformation of the form 
$$[x : y : z] \mapsto [P(x,y,z) : Q(x,y,z) : R(x,y,z)],$$
\noindent where $P,Q,R$ are three homogeneous polynomials of the same degree, and which admits an inverse of the same form. It was proved by Noether that $\bir$ is generated by  the standard quadratic involution
$$\sigma : [x : y : z] \mapsto [\frac{1}{x}: \frac{1}{y}:\frac{1}{z}]=[yz : xz : xy],$$

\noindent and the group $\PGL$ of automorphisms of $\mathbb{P}^{2}$. This group was studied a lot by algebraic geometers, and we refer the reader to \cite{cantat} for some references. Recently, Cantat and Deserti started the study of infinite finitely generated subgroups of $\bir$, using ideas from geometric group theory, see \cite{cantat,cantatlamy,deserti2,deserti1} as well as \cite{favre}. In particular, using ideas of Zariski and Manin, Cantat recently proved that the group $\bir$ can be embedded into the isometry group of an infinite dimensional hyperbolic space $\mathbb{H}_{\mathbb{P}^{2}}$, called the {\it Picard-Manin space}. We will briefly recall its construction in section~\ref{section-cremo}. Using this embedding he proved the following result: {\it any homomorphism with infinite image from a discrete Kazhdan group into the Cremona group $\bir$ is conjugated to a homomorphism into $\PGL$}. In particular, this result applies to any lattice $\Lambda$ in a connected simple Lie group with property ${\rm T}$. This left open the problem of classifying homomorphisms from lattices in the groups ${\rm SO}(n,1)$ and $\SU$ into the Cremona group. Note that this discussion is motivated by the fact that (for some values of $n$) there exist injective homomorphisms from lattices in ${\rm SO}(n,1)$ to the Cremona group (see \cite{cantatlamy,dolzh} and the references there).

\begin{rem} We will see that the space $\mathbb{H}_{\mathbb{P}^{2}}$ is not separable; however any finitely generated group acting on it preserves a closed, separable, totally geodesic subspace. 
\end{rem}

In the following, we say that a homomorphism from a finitely generated group into the group $\bir$ is non-elementary if the associated action on $\mathbb{H}_{\mathbb{P}^{2}}$ is non-elementary. In the same way, we say that a subgroup of the Cremona group is elliptic if it fixes a point in the space $\mathbb{H}_{\mathbb{P}^{2}}$. Using our previous theorem, and the fact that the action of $\bir$ on $\mathbb{H}_{\mathbb{P}^{2}}$ is, in some sense, discrete, we can prove:

\begin{main}\label{cremo} Let $\Gamma$ be the fundamental group of a compact K\"{a}hler manifold $X$. Let $$\rho : \Gamma \to \bir$$ be a non-elementary homomorphism. Then there exists a fibration $X\to \Sigma$ onto a hyperbolic $2$-orbifold such that the kernel $H$ of the map $\Gamma \to \pi_{1}^{orb}(\Sigma)$ has the following property. The group $\rho(H)$ fixes pointwise a closed totally geodesic subspace of $\mathbb{H}_{\mathbb{P}^{2}}$ of dimension at least $2$. \end{main}

Observe that if the group $\rho(H)$ appearing in the theorem above is not trivial, it provides an example of an elliptic subgroup of the Cremona group for which the action of the normalizer is non-elementary (since the action of $\Gamma$ itself is non-elementary). To study further homomorphisms from K\"{a}hler groups to the Cremona group, one is thus led to study elliptic subgroups of the Cremona group whose normalizer is large enough. Such an elliptic subgroup can be infinite as the following example shows. 

\begin{exam} In $\mathbb{P}^{2}$, with the homogeneous coordinates $[x:y:z]$, consider the complement of the union of the three lines $\{x=0\}$, $\{y=0\}$ and $\{z=0\}$. This is an open set $O$ isomorphic to $\mathbb{C}^{\ast}\times \C^{\ast}$. Thinking of $O$ as a group, we get an action of the group $\C^{\ast}\times \C^{\ast}$ on $O$ by translation. On the other hand, the group ${\rm GL}(2,\mathbb{Z})$ acts on $O$ by monomial transformations: the matrix 
$$\left( 
\begin{array}{cc}
a & b\\
c & d\\
\end{array}\right)$$

\noindent acts by the transformation $(x,y)\mapsto (x^{a}y^{b},x^{c}y^{d})$. One obtains in this way an injective homomorphism of the group $\left(\C^{\ast}\times \C^{\ast}\right)\rtimes {\rm GL}(2,\mathbb{Z})$ into the Cremona group. We will denote by $G_{toric}$ its image. One can see that the group $\mathbb{C}^{\ast}\times \C^{\ast}$ fixes pointwise an infinite dimensional totally geodesic subspace of $\mathbb{H}_{\mathbb{P}^{2}}$. We will explain this in section~\ref{section-cremo}, after constructing the space $\mathbb{H}_{\mathbb{P}^{2}}$.
\end{exam}

This is essentially the only example of an elliptic subgroup of $\bir$ with large normalizer as shown by the following theorem proved by Serge Cantat in the appendix: 
\begin{main}\label{serge-elliptic}
Let $N$ be a subgroup of the Cremona group $\bir$. Assume that there exists a short exact
sequence 
\[
1 \to A \to N \to B\to 1
\]
where $N$ contains at least one hyperbolic element,  and $A$ is infinite elliptic. Then $N$ is conjugate to a subgroup of the group $G_{toric}$ of automorphisms of $(\C^*)^2$.
\end{main}

Note that, as in the finite dimensional case, isometries of infinite dimensional hyperbolic spaces fall into three types: elliptic, hyperbolic and parabolic; see~\cite{bim}. The word {\it hyperbolic} in the theorem above refers to the type of the action of a birational map on the space $\mathbb{H}_{\mathbb{P}^{2}}$.

 Combining Theorem~\ref{cremo} with Cantat's theorem, we obtain a more precise description of non-elementary homomorphisms  from K\"{a}hler groups to the Cremona group:

\begin{cor}\label{nonelecomplet} Let $\rho : \Gamma \to \bir$ be a non-elementary homomorphism. One of the following two cases occurs:
\begin{enumerate}
\item The homomorphism $\rho$ is conjugated to a homomorphism with values in the group $G_{toric}$.
\item After passing to a finite index subgroup of $\Gamma$ (equivalently, after taking a finite cover of $X$), $\rho$ factors through a fibration onto a hyperbolic $2$-orbifold. 
\end{enumerate}
\end{cor}
{\it Proof.} We keep the notation of Theorem~\ref{cremo}. If the group $\rho(H)$ is infinite, Cantat's theorem implies that $\rho(\Gamma)$ is conjugated to a subgroup of the group $G_{toric}$. If $\rho(H)$ is finite, we argue as follows. We consider the homomorphism $\Gamma \to {\rm Aut}(\rho(H))$ given by the action of $\Gamma$ on $\rho(H)$ by conjugation. Let $\Gamma_{1}$ be its kernel, $H_{1}$ be the intersection of $\Gamma_{1}$ with $H$ and $\Lambda_{1}$ be the image of $\Gamma_{1}$ in $\pi_{1}^{orb}(\Sigma)$. Up to replacing $\Gamma_{1}$ by another finite index subgroup we also assume that $\Lambda_{1}$ is torsion free, hence a surface group. 

Let $V=\rho(H_{1})$; this is a finite abelian subgroup of the Cremona group, and the {\it central} extension $$1 \to V \to \rho(\Gamma_{1})\to \rho(\Gamma_{1})/V\to 1$$
\noindent determines a class in $H^{2}(\rho(\Gamma_{1})/V,V)$. Let $e$ be the pull-back of this class in $H^{2}(\Lambda_{1},V)$.  The class $e$ becomes trivial on a finite index subgroup $\Lambda_{2}$ of $\Lambda_{1}$. Let $\Gamma_{2}$ be the inverse image of $\Lambda_{2}$ in $\Gamma_{1}$. The fact that the class $e$ is trivial on $\Lambda_{2}$ says that the homomorphism $\rho : \Gamma_{2} \to \bir$ factors through a homomorphism to the direct product $\Lambda_{2} \times V$. Hence, up to taking once again a finite index subgroup, it factors through the projection on $\Lambda_{2}$.\hfill $\Box$

The following corollary answers partially a question raised by Cantat~\cite{cantat}:

\begin{cor}\label{sc} Let $\Gamma_{1}$ be a cocompact lattice in the group $\SU$ with $n\ge 2$. If $$\rho : \Gamma_{1}\to \bir$$ is an injective homomorphism, then one of the following two possibilities holds:
\begin{enumerate}
\item The group $\rho(\Gamma_{1})$ fixes a point in the Picard-Manin space $\mathbb{H}_{\mathbb{P}^{2}}$.
\item The group $\rho(\Gamma_{1})$ fixes a unique point in the boundary of the Picard-Manin space $\mathbb{H}_{\mathbb{P}^{2}}$.
\end{enumerate}
\end{cor}

The proof of Corollary~\ref{sc} will be given in Section~\ref{section-cremo}.

A natural problem left open by this work is the study of homomorphisms from K\"{a}hler groups to the group $G_{toric}$ as well as homomorphisms into the Cremona group for which the associated action on $\mathbb{H}_{\mathbb{P}^{2}}$ has a unique fixed point at infinity (see~\cite{cantat} for an interpretation of these morphisms in terms of algebraic geometry). It might also be possible to extend these results to noncocompact lattices in $\SU$ (or more generally fundamental groups of quasiprojective varieties) using techniques similar to the ones used in \cite{corsimp} and \cite{komau}, but we have not tried to establish such results.

\subsection{About the proofs}

 We now briefly discuss the (classical) strategy of the proof of Theorem~\ref{infinite}. We denote by $\widetilde{X}$ the universal cover of $X$. We use the fact that there is an equivariant harmonic map $f$ from $\widetilde{X}$ to $\HHI$, this is due to Korevaar and Schoen \cite{korschoen}. It is well-known \cite{sampson2} (at least in the finite-dimensional case) that this map has real rank at most $2$. This is the main reason why actions of K\"{a}hler groups on real hyperbolic spaces are so constrained. We can then prove a {\it factorization theorem} in the spirit of many factorization theorems in the literature (\cite{ct1,gs,jostyau,jostzuo}...), i.e. we prove that $f$ can be written as the composition of a holomorphic map from $\widetilde{X}$ to the unit disc $\Delta \subset \mathbb{C}$ and a harmonic map from $\Delta$ to $\HHI$. A similar result appears in \cite{jostyau}, however we give a detailed account of the proof here. This decomposition gives rise to a homomorphism $\theta$ from $\Gamma$ to $\PSLR$. Our main observation, inspired by the reading of \cite{gs} and which is essentially the only new remark in the proof, is that if the group $\theta(\Gamma) \subset \PSLR$ is not discrete, the natural action of $\theta(\Gamma)$ on $\HHI$ extends to a continuous action of $\PSLR$. This allows to prove easily Theorem~\ref{infinite}.

\begin{rem}\label{symspace} The same arguments show that a Zariski dense representation of a K\"{a}hler group $\Gamma$ in a noncompact simple Lie group $G$ (with trivial center and not isomorphic to $\PSLR$) will factor through a fibration onto a hyperbolic $2$-orbifold as soon as the associated harmonic map $f$ to the symmetric space of $G$ has the following property: the complex linear part $df^{1,0}$ of the differential of $f$ has complex rank $1$ (see the end of section~\ref{conclusionoftheproof}).  
\end{rem}

The text is organized as follows. Section~\ref{sl2} contains the description of exotic actions of $\PSLR$ on $\HHI$ and the proof of Theorem~\ref{minim}. In section~\ref{section-finite}, we recall a few facts concerning harmonic maps, for the reader interested in the application to the Cremona group and not familiar with these techniques. We also check that classical results (such as the Bochner-Siu-Sampson formula) still apply when the target of the harmonic map is an infinite dimensional hyperbolic space. In section~\ref{section-infinite}, we prove Theorem~\ref{infinite}. The application to the Cremona group is described in section~\ref{section-cremo}. In section~\ref{section-dbar}, we check that the classical Koszul-Malgrange theorem concerning integrability of certain almost complex structures on complex vector bundles can be adapted to the infinite-dimensional setting. This is needed in section~\ref{section-infinite}. Finally, section~\ref{analyse} contains the proof of a technical proposition needed in section~\ref{section-finite}. 

\noindent {\bf Acknowledgements.} We would like to thank \'Etienne Ghys, who first told us about exotic actions of $\PSLR$ on $\HHI$, Serge Cantat for his explanations concerning the Cremona group and for his comments on the text, and Nicolas Monod for several conversations and for pointing out to us an inaccuracy in the first version of this text. We also want to thank Olivier Biquard, Dominique Cerveau, Benson Farb, Daniel Panazzolo and Robert Stanton for several useful conversations, as well as the referee for his comments on the text.

%%%%%%%%%%%%%%%%%%%%%%%%%%%%%%%%%%%%%%%%%%%%%%%%%%%%%%%%%%%%%%%%%%%%%%%%%%%%%%%%%%%%%%%
%%%%%%%%%%%%%%%%%%%%%%%%%%%%%%%%%%%%%%%%%%%%%%%%%%%%%%%%%%%%%%%%%%%%%%%%%%%%%%%%%%%%%%%

\section{Exotic actions of $\PSLR$}\label{sl2}

In this section we describe some infinite dimensional linear representations of $\PSLR$ which lead to a family of exotic actions of $\PSLR$ on the infinite dimensional real hyperbolic space. This will lead to the proof of Theorem~\ref{minim}. 

\subsection{A glimpse at the representation theory of $\PSLR$}

Let $\Delta$ be the unit disc in $\C$. We will consider the usual action of $\PSLR$ on $\Delta$ obtained by conjugating the homographic action of $\PSLR$ on the upper half plane by the map $z\mapsto \frac{z-i}{z+i}$. All along this section, we will denote by $\mathscr{H}_{\mathbb{C}}$ the Hilbert space of square integrable functions on ${\rm S}^{1}=\partial \Delta$, endowed with the angular measure $d\theta$ and for $n\in \Z$ by $e_{n}$ the function $z\mapsto z^{n}$ on the circle. For each real number $s$, we define a representation $\pi_{s}$ of $\PSLR$ on $\mathscr{H}_{\mathbb{C}}$ by:
$$\pi_{s}(g)\cdot f=Jac(g^{-1})^{\frac{1}{2}+s}f\circ g^{-1}\;\;\;\;\; (g\in \PSLR, f\in \mathscr{H}_{\mathbb{C}}).$$

\noindent Here $Jac(g)$ is the jacobian of an element $g$ with respect to the measure $d\theta$ on the circle. For each $g\in \PSLR$ the operator $\pi_{s}(g) : \mathscr{H}_{\mathbb{C}}\to \mathscr{H}_{\mathbb{C}}$ is a bounded linear operator (although not unitary in general). This defines a continuous representation of $\PSLR$ on $\mathscr{H}_{\mathbb{C}}$. It is well-known that the representation $\pi_{s}$ is irreducible if and only if $s\notin \frac{1}{2}+\mathbb{Z}$ (see \cite{knapp} or \cite{wallach} for instance).

\begin{rem} Note that the formula above still defines a representation when $s$ is a complex number. When $s$ is purely imaginary the representation is unitary and one obtains the so-called {\it principal series} of unitary representations of $\PSLR$ (although our description of these representations might differ from the usual one by a translation in the parameter~$s$). 
\end{rem}

If $f_{1},f_{2}\in \mathscr{H}_{\mathbb{C}}$, we will write $\left( f_{1},f_{2}\right):=\int_{{\rm S}^{1}}f_{1}\overline{f_{2}}d\theta$. From the definition of the jacobian, we obtain:

$$\left( \pi_{s}(g)(f_{1}),\pi_{-s}(g)(f_{2})\right)=\left( f_{1},f_{2}\right)\;\;\;\;\; (f_{1},f_{2}\in \mathscr{H}_{\mathbb{C}}, g\in \PSLR).$$

\noindent Hence, if one can construct a continuous operator $A_{s}$ from $\mathscr{H}_{\mathbb{C}}$ to itself which intertwines the representations $\pi_{s}$ and $\pi_{-s}$ (i.e. satisfies $\pi_{-s}(g)\circ A_{s}=A_{s}\circ \pi_{s}(g)$), this will imply that there is a hermitian pairing on $\mathscr{H}_{\mathbb{C}}$, invariant by the representation $\pi_{s}$. It turns out that such an operator exists (for the right values of $s$) and has been much studied in the litterature on representation theory, see for instance \cite{knapp,sally1,sally2}, as well as \cite{jonhsonwallach,knapp0} for a study of intertwining operators for representations of more general simple Lie groups than $\PSLR$. Since this material might not be so well-known outside the representation theory community, we will include a brief exposition of the case of interest to us. Namely, we will prove the following proposition:

\begin{prop}\label{entrelace} Let $s>0$ and let $A_{s} : \mathscr{H}_{\mathbb{C}} \to \mathscr{H}_{\mathbb{C}}$ be the operator which maps $e_{n}$ to $\lambda_{n}e_{n}$; where $\lambda_{n}$ is given by $\lambda_{0}=1$ and for $n\neq 0$:
$$\lambda_{n}=\prod_{i=0}^{\vert n\vert -1}\frac{i+\frac{1}{2}-s}{i+\frac{1}{2}+s}.$$
\noindent Then $A_{s}$ defines a bounded linear operator from $\mathscr{H}_{\mathbb{C}}$ to $\mathscr{H}_{\mathbb{C}}$ which satisfies $\pi_{-s}(g)\circ A_{s}=A_{s}\circ \pi_{s}(g)$ for all $g\in \PSLR$. 
\end{prop}

The reader will notice that the formula for $A_{s}$ given above does not define a bounded linear operator on $\mathscr{H}_{\mathbb{C}}$ if $s<0$. Indeed if $A_{s}$ was bounded for some negative $s$ not in $\frac{1}{2}+\mathbb{Z}$, there would exists a constant $C>0$ such that $$\left| \prod_{i=0}^{\vert n\vert -1}\frac{i+\frac{1}{2}-s}{i+\frac{1}{2}+s}\right| \le C.$$
\noindent This would imply the convergence of the series $$\sum_{i\ge 0} \log(1+\frac{-2s}{i+\frac{1}{2}+s}),$$ a contradiction. This implies that for $s>0$, $s\notin \frac{1}{2}+\mathbb{Z}$, $A_{s}$ has a dense, {\it non-closed} image: if $A_{s}$ was surjective, it would be an isomorphism by the open mapping theorem. Its inverse $A_{s}^{-1}=A_{-s}$ would then be bounded.

Once the above proposition is known we can define for $f_{1},f_{2}\in \mathscr{H}_{\mathbb{C}}$:
$$\langle f_{1},f_{2}\rangle_{s}=\int_{S^{1}}f_{1}\overline{A_{s}(f_{2})}d\theta.$$
\noindent The bilinear form $\langle \cdot , \cdot \rangle_{s}$ on $\mathscr{H}_{\mathbb{C}}$ is invariant by the representation $\pi_{s}$. If $s\in (0,\frac{1}{2})$, it is positive definite. This implies that the representations $(\pi_{s})_{0<s<\frac{1}{2}}$ are unitary (although the invariant inner product is not the standard one on $\mathscr{H}_{\mathbb{C}}$). This is the so-called {\it complementary series} which arises in the classification of irreducible unitary representations of $\PSLR$. For $s>\frac{1}{2}$ the pairing $\langle \cdot , \cdot \rangle_{s}$ is not positive definite anymore, but this does not make it less interesting! If $s\in (p-\frac{1}{2}, p+\frac{1}{2})$ where $p\ge 1$ is an integer, it is easily checked that the bilinear form $\langle\cdot , \cdot \rangle_{s}$ has index $p$ (here {\it index} should be understood as the dimension of a maximal isotropic subspace):

\begin{itemize}
\item If $s\in (p-\frac{1}{2}, p+\frac{1}{2})$ and $p$ is odd, $p$ of the coefficients $(\lambda_{n})_{n\in \mathbb{Z}}$ are positive, all the others are negative. 
\item If $s\in (p-\frac{1}{2}, p+\frac{1}{2})$ and $p$ is even, $p$ of the coefficients $(\lambda_{n})_{n\in \mathbb{Z}}$ are negative, all the others are positive.
\end{itemize} 

In particular, for each $s\in (\frac{1}{2},\frac{3}{2})$, we can look at the restriction of the representation $\pi_{s}$ to the subspace $\mathscr{H}_{\mathbb{R}}\subset \mathscr{H}_{\mathbb{C}}$ of real-valued functions. However we are not yet in the situation described in the introduction to construct the space $\HHI$. Indeed, the space of functions with zero mean, endowed with the scalar product $-\langle \cdot , \cdot \rangle_{s}$ is not complete. This follows from the fact that the sequence $(\lambda_{n})_{n\ge 0}$ which appears in Proposition~\ref{entrelace} above tends to zero as $n$ goes to $+\infty$. 

If $f\in \mathscr{H}_{\mathbb{R}}$, we will denote by $f_{0}$ the function 
$$f-\frac{1}{2\pi}\int_{S^{1}}fd\theta,$$
\noindent and by $m(f)$ the number $\int_{S^{1}}fd\theta$. Let $\overline{\mathscr{H}_{\mathbb{R}}}$ be the completion of the space $\mathscr{H}_{\mathbb{R}}$ for the scalar product associated to the norm
$$\vert \vert f\vert \vert =\sqrt{\left( \int_{S^{1}}fd\theta \right)^{2}-\langle f_{0},f_{0}\rangle_{s}}.$$
\noindent The bilinear form $\langle \cdot , \cdot \rangle_{s}$ extends to $\overline{\mathscr{H}_{\mathbb{R}}}$ and its extension is still denoted by the same symbol. We then have the:

\begin{prop} The action of $\PSLR$ on $\mathscr{H}_{\mathbb{R}}$ extends to a continuous, irreducible, linear representation $\overline{\pi_{s}}$ on $\overline{\mathscr{H}_{\mathbb{R}}}$. The action of $\PSLR$ on $\overline{\mathscr{H}_{\mathbb{R}}}$ preserves the extended quadratic form $\langle \cdot , \cdot \rangle_{s}$. 
\end{prop}
\noindent {\it Proof.} Let $g$ be in $\PSLR$. We first look for a constant $C$ such that $$\vert \vert \pi_{s}(g)(f)\vert \vert \le C \vert \vert f\vert \vert,$$\noindent for $f\in \mathscr{H}_{\mathbb{R}}$. Write $\pi_{s}(g)(f_{0})=c+u_{0}$ where $c$ is constant and $u_{0}$ is a function with zero mean. We have:
$$\vert c\vert =\vert \langle \pi_{s}(g)(f_{0}),(2\pi)^{-1}\rangle_{s}\vert =\vert \langle f_{0},\pi_{s}(g^{-1})((2\pi)^{-1})\rangle_{s}\vert\le \vert \vert f_{0}\vert \vert \cdot C_{1} \le \vert \vert f\vert \vert \cdot C_{1},$$
\noindent where $C_{1}=\vert \vert (\pi_{s}(g^{-1})((2\pi)^{-1}))_{0}\vert \vert$. Since $\pi_{s}(g)$ preserves $\langle \cdot , \cdot \rangle_{s}$, we have:
$$-\langle u_{0},u_{0}\rangle_{s}=c^{2}-\langle f_{0},f_{0}\rangle_{s}\le c^{2}+\vert \vert f\vert \vert^{2}\le (1+C_{1}^{2})\vert \vert f\vert \vert^{2}.$$ \noindent Finally we have:
$$\vert \vert \pi_{s}(g)(f)\vert \vert \le \vert \vert \pi_{s}(g)(m(f))\vert \vert +\vert \vert c\vert \vert +\vert \vert u_{0} \vert \vert \le C \vert \vert f\vert \vert,$$

\noindent for $C=\vert \vert \pi_{s}(g)(1)\vert \vert+C_{1}+\sqrt{1+C_{1}^{2}}$. This implies that $\pi_{s}(g)$ extends to a bounded operator $$\overline{\pi_{s}}(g) : \overline{\mathscr{H}_{\mathbb{R}}}\to \overline{\mathscr{H}_{\mathbb{R}}}.$$ \noindent One deduces easily from this that $\overline{\pi_{s}}$ is continuous and preserves the (extension of the) form $\langle \cdot , \cdot \rangle_{s}$ on $\overline{\mathscr{H}_{\mathbb{R}}}$. As for the irreducibility of $\overline{\pi_{s}}$, we argue as follows. Let $V\subset \overline{\mathscr{H}_{\mathbb{R}}}$ be a nonzero closed $\PSLR$-invariant subspace. It is in particular invariant under the group ${\rm SO}(2)$, which preserves the norm $\vert \vert \cdot \vert \vert$. Hence $V$ is a direct sum of $1$-dimensional subspaces invariant by ${\rm SO}(2)$. From this, one sees that $V$ must contain a nonzero vector in $\mathscr{H}_{\mathbb{R}}$. The irreducibility of $\pi_{s}$ then implies that $V$ contains $\mathscr{H}_{\mathbb{R}}$. Hence $V=\overline{\mathscr{H}_{\mathbb{R}}}$ since it is closed.\hfill $\Box$

As a consequence of the previous proposition, for each $s\in (\frac{1}{2},\frac{3}{2})$, the representation $\overline{\pi_{s}}$ defines a non-elementary action of $\PSLR$ on the hyperbolic space $\HHI$, with no nontrivial closed invariant totally geodesic subspace. When the parameter $s$ if greater than $\frac{3}{2}$, one obtains similarly actions of $\PSLR$ on the infinite dimensional symmetric spaces associated to the (naturally defined) groups ${\rm O}(p,\infty)$ (see \cite{gromov1}, section 6, for a description of these spaces). 

\begin{rem} The article \cite{jonhsonwallach} contains similar results for all classical simple Lie groups of rank $1$ (namely ${\rm SO}(n,1)$, ${\rm SU}(n,1)$ and ${\rm Sp}(n,1)$). For ${\rm SO}(n,1)$ ($n\ge 3$) one obtains invariant quadratic forms of finite index on certain representations, the index depending on the dimensions of the space of spherical harmonics of degree $k$. Note that the index $1$ can be obtained. For ${\rm SU}(n,1)$ and ${\rm Sp}(n,1)$ however, the quadratic forms that one obtains are of infinite index (when they are not unitary), see the formulas on page 156 of \cite{jonhsonwallach}.  
\end{rem}

{\it Proof of Proposition~\ref{entrelace}.} To prove the proposition, we use a different model for the representation $\pi_{s}$ that we now describe. Let $\xi =1\in \partial \Delta$ be the unique point fixed by the affine group of matrices of the form:
$$\left( \begin{array}{cc}\ast & \ast'\\ 0 & \ast^{-1}\\ \end{array}\right).$$ 
\noindent We will use the fact that any element $g\in {\rm SL}_{2}(\mathbb{R})$ can be uniquely written as: 
$$g=ka_{\lambda}n_{t},$$

\noindent where $k\in {\rm SO}(2)$, $a_{\lambda}=\left( \begin{array}{cc}\lambda & 0\\ 0 & \lambda^{-1}\\ \end{array}\right)$ and $n_{t}=\left( \begin{array}{cc}1 & t \\ 0 & 1\\ \end{array}\right)$. Let $X^{s}$ be the space of measurable functions $F : {\rm SL}_{2}(\mathbb{R})\to \C$ which satisfy the following properties: 
\begin{itemize}
\item $F(-g)=F(g)$,
\item the restriction of $F$ to ${\rm SO}(2)$ is square integrable,
\item $F(ga_{\lambda}n_{t})=F(g)\lambda^{-1-2s}$. 
\end{itemize}
In the following, we identify ${\rm SO}(2)/\pm {\rm Id}$ to $\partial \Delta$ via the map $k\mapsto k(\xi)$. Now, if $F\in X^{s}$, its restriction to ${\rm SO}(2)$ defines an element of $\mathscr{H}_{\mathbb{C}}$ and conversely, if $f\in \mathscr{H}_{\mathbb{C}}$, one can extend it to a function $F \in X^{s}$ by declaring that $F(ka_{\lambda}n_{t})=f(k(\xi))\lambda^{-1-2s}$. One verifies that under the previous identification of $X^{s}$ with $\mathscr{H}_{\mathbb{C}}$ the action of $\pi_{s}(g)$ becomes simply the precomposition with the left translation by $g^{-1}$. Consider now the operator $F \mapsto LF$ defined by:
$$LF(x)=\int_{\R}F(xw\underline{n}(u))du,$$

\noindent where $\underline{n}_{u}=\left( \begin{array}{cc} 1 & 0\\ u & 1\\ \end{array}\right)$ and $w=\left( \begin{array}{cc} 0 & -1\\ 1 & 0\\ \end{array}\right)$. It is easy to check that if $F\in X^{s}$ is continuous, the previous integral is well-defined and defines a continuous function on ${\rm SL}_{2}(\mathbb{R})$. We will check that $LF\in X^{-s}$. Observe that one has the following identity:
$$xa_{\lambda}n_{t}w\underline{n}(u)=xw\underline{n}(\lambda^{2}(u-t)) a_{\lambda^{-1}};$$

\noindent hence:
$$F(xa_{\lambda}n_{t}w\underline{n}(u))=F(xw\underline{n}(\lambda^{2}(u-t)))\lambda^{1+2s},$$
\noindent and a change of variable gives: $LF(xa_{\lambda}n_{t})=LF(x)\lambda^{-(1-2s)}$. To finish the proof of the proposition, we thus have to check that the operator $L$ extends continuously to all of $X^{s}$ and coincides (up to a scalar factor) with the operator $A_{s}$ defined by the formula in the statement of the proposition. Note that by its very definition, the operator $L$ intertwines the actions of ${\rm SL}(2,\R)$ on $X^{s}$ and $X^{-s}$. Let $k_{\theta}\in {\rm SO}(2)$ be the rotation of angle $\theta$. A simple calculation shows that:

$$w\underline{n}(u)=k_{\theta_{u}}a_{\sqrt{1+u^{2}}}n_{\frac{u}{1+u^{2}}}$$

\noindent where the angle $\theta_{u}$ is defined by the relation $e^{i\theta_{u}}=\frac{i-u}{\sqrt{1+u^{2}}}$. From this we deduce that $e_{k}(k_{\theta}w\underline{n}(u))=e_{k}(k_{\theta})(\frac{i-u}{\sqrt{1+u^{2}}})^{2k}(1+u^{2})^{-\frac{1}{2}-s}$. This implies that $Le_{k}=\lambda_{k}e_{k}$, where $\lambda_{k}$ is given by:
$$\lambda_{k}=\int_{\R}    \left( \frac{i-u}{\sqrt{1+u^{2}}}\right)^{2k}\frac{du}{(1+u^{2})^{\frac{1}{2}+s}}=\int_{\R}\frac{(i-u)^{2k}du}{(1+u^{2})^{k+\frac{1}{2}+s}}.$$ All $\lambda_{k}$'s having modulus bounded by $\int_{\R}\frac{du}{(1+u^{2})^{\frac{1}{2}+s}}$, $L$ extends continuously to $X^{s}$. Since $\lambda_{k}=\lambda_{-k}$ it is now enough to check the equality $\lambda_{k+1}=\frac{k+\frac{1}{2}-s}{k+\frac{1}{2}+s}\lambda_{k}$ for $k\ge 0$. Up to dividing by $\lambda_{0}$, this will show that the operator $L$ has the form given in the statement of the proposition. We check this equality in the next lemma.\hfill $\Box$

\begin{lemma} For any positive integer $k$, one has: $\lambda_{k+1}=\frac{k+\frac{1}{2}-s}{k+\frac{1}{2}+s}\lambda_{k}$. 
\end{lemma}
{\it Proof.} We start from the equality $\lambda_{k}=\int_{\R}\frac{(i-u)^{2k}}{(1+u^{2})^{k+\frac{1}{2}+s}}du$ and integrate by part (integrating $(i-u)^{2k}$) to obtain:
$$\begin{array}{ccc}
\lambda_{k} & = & \frac{2(k+\frac{1}{2}+s)}{2k+1}\left( \lambda_{k+1}-\int_{\R}\frac{(i+u)(i-u)^{2k+1}}{(1+u^{2})^{k+\frac{3}{2}+s}}du+\int_{\R}\frac{u(i-u)^{2k+1}}{(1+u^{2})^{k+\frac{3}{2}+s}}du\right). \\
 \end{array}$$

\noindent The second term in the parenthesis above equals $\lambda_{k}$; the third term is proportional to $\lambda_{k}$ (by an integration by part which is inverse to the one we just performed). Hence we get:
$$\lambda_{k}=\frac{2(k+\frac{1}{2}+s)}{2k+1}\left(\lambda_{k+1}+\lambda_{k}-\frac{2k+1}{2(k+\frac{1}{2}+s)}\lambda_{k}\right)$$
\noindent The proof is now finished by a simple calculation.\hfill $\Box$

\subsection{Some properties of actions of $\PSLR$ on $\HHI$}

In this section we prove a proposition which describes how the family $(\overline{\pi_{s}})_{s\in (\frac{1}{2},\frac{3}{2})}$ of actions of $\PSLR$ on $\HHI$ varies with $s$. This will complete the proof of Theorem~\ref{minim} (setting $t=s-\frac{1}{2}$, one recovers the notation of the theorem).

We denote here by $\HHI_{s}$ the hyperbolic space constructed from the pair $(\overline{\mathscr{H}_{\mathbb{R}}},\langle \cdot , \cdot \rangle_{s})$, where $\frac{1}{2}<s<\frac{3}{2}$ (all these spaces are isometric but we will keep this convenient notation). Observe that the group ${\rm SO}(2)$ has a unique fixed point in $\HHI_{s}$, hence there exists a unique $\overline{\pi_{s}}$-equivariant map $f_{s} : \Delta \to \HHI_{s}$. This map is automatically harmonic: indeed if the tension field of $f_{s}$ was nonzero at a point $\xi \in \Delta$ the stabilizer of $\xi$ in $\PSLR$ would fix the geodesic going through $f_{s}(\xi)$ and tangent to the tension field there, a contradiction. For the definition of the tension field of a map, or of its harmonicity, we refer the reader to paragraph~\ref{classicalcase}.

\begin{prop} If $g_{hyp}$ denotes the hyperbolic metric with curvature $-1$ both on $\HHI_{s}$ and on $\Delta$, we have: $f_{s}^{\ast}g_{hyp}=c_{s}g_{hyp}$ where:
$$c_{s}=\frac{(1+2s)(s-\frac{1}{2})}{4}.$$
\end{prop}

As a consequence of this proposition, the image of $f_{s}$ is a minimal surface in $\HHI_{s}$ with curvature $-c_{s}^{-1}$ (which goes to $-\infty$ or $-1$ as $s$ goes to $\frac{1}{2}$ or $\frac{3}{2}$ respectively).  

{\it Proof.} Since there is only one $\PSLR$-invariant metric on the unit disc, up to a scalar, we can write: $f_{s}^{\ast}g_{hyp}=c_{s}g_{hyp}$, for some constant $c_{s}$. We want to compute the value of $c_{s}$. It will be convenient to identify $\PSLR$ with ${\rm SU}(1,1)$ (up to a $2$-sheeted cover). So let $g_{t} : \Delta \to \Delta$ be the isometry induced by the matrix
$$\left( \begin{array}{cc}
{\rm ch}(t) & {\rm sh}(t)\\
{\rm sh}(t) & {\rm ch}(t)\\
\end{array}\right).$$
\noindent We will also denote by $o$ the origin of the disc. We will prove the following estimate:
$$d_{\HHI_{s}}(f_{s}(g_{t}(o)),f_{s}(o))=\sqrt{(1+2s)(s-\frac{1}{2})}t+o(t).\;\;\;\; (\ast)$$
This allows to conclude easily. Denote by $v$ the derivative of the curve $g_{t}\cdot o$ (which is nothing else than $1\in \mathbb{C}$). Since the hyperbolic metric on the disc is $\frac{4dzd\overline{z}}{(1-\vert z \vert^{2})^{2}}$ we have $g_{hyp}(v,v)=4$. On the other hand the equality above says that $g_{hyp}(Df_{s}(v),Df_{s}(v))=(1+2s)(s-\frac{1}{2})$. Hence $c_{s}=\frac{(1+2s)(s-\frac{1}{2})}{4}$.

We now prove the equality $(\ast)$. Observe that the constant function equal to $\frac{1}{\sqrt{2\pi}}$ is the unique fixed point of the group ${\rm SO}(2)$ in $\HHI_{s}$, hence $f_{s}(o)=\frac{1}{\sqrt{2\pi}}$ and $f_{s}(g_{t}(o))$ is the function sending a point $e^{i\theta}$ to $(2\pi)^{-\frac{1}{2}}\vert {\rm sh}(t)e^{i(\theta+\pi)}+{\rm ch}(t)\vert^{-1-2s}$. Define:
$$u(s,t):=\langle \pi_{s}(g_{t})(1),1\rangle_{s}=\frac{1}{2\pi}\int_{0}^{2\pi}\frac{d\theta}{\vert {\rm sh}(t)e^{i\theta}+{\rm ch}(t)\vert^{1+2s}}.$$

\noindent A calculation shows that:
$$\frac{1}{\vert {\rm sh}(t)e^{i\theta}+{\rm ch}(t)\vert^{1+2s}}=1-(\frac{1}{2}+s)2t\cos (\theta)+((1+2s)(\frac{3}{2}+s)\cos^{2}(\theta)-(1+2s))t^{2}+o(t^{2})$$
\noindent and then that:
$$u(s,t)=1+(1+2s)(\frac{s}{2}-\frac{1}{4})t^{2}+o(t^{2}).$$

\noindent Remembering that $u(s,t)={\rm ch}(d(\pi_{s}(g_{t})\cdot 1,1))$, one easily deduces the formula~$(\ast)$.\hfill $\Box$

%%%%%%%%%%%%%%%%%%%%%%%%%%%%%%%%%%%%%%%%%%%%%%%%%%%%%%%%%%%%%%%%%%%%%%%%%%%%%%%%%%%%%%%%%%%%%%%%%%%%%%%%%%%%%%%%%%%%%%%%%%%%%%%%%%%%%%%%%%%%%%%%%%%%%%%%%%%%%%%%%%%%%%%%%%%%%%%%

\section{Preliminaries on harmonic maps}\label{section-finite}

\subsection{The classical case}\label{classicalcase}

We start by recalling a few properties of {\it harmonic maps}. Some references for a more detailed presentation of the material in this section are \cite{period,gs,korschoen0,nishi}. 

A map $f : M \to N$ between two Riemannian manifolds is harmonic if it is a critical point of the energy functional $E(f):=\int_{M}e(f)(x)dx$, where $dx$ is the Riemannian measure on $M$ and where $e(f)(x)=\frac{1}{2}\vert \vert df_{x}\vert \vert^{2}$ is the energy density of $f$. The {\it tension field} of a map $f : M \to N$ is the vector field along the image of $f$ defined as follows. Let $\nabla$ be the connection on the vector bundle ${\rm T}^{\ast}M\otimes f^{\ast}{\rm T}N$ constructed from the Levi-Civita connections of $M$ and $N$. The differential $df$ of $f$ can be thought of as a section of the bundle ${\rm T}^{\ast}M\otimes f^{\ast}{\rm T}N$, and the tension field $\tau (f)$ of $f$ is the trace of the $2$-tensor $\nabla df$:
$$\tau(f)=\sum_{i=1}^{m} \nabla df(e_{i},e_{i}),$$
\noindent where $(e_{i})_{1\le i\le m}$ is a local orthonormal frame of ${\rm T}M$. One can show that the harmonicity of $f$ is equivalent to the vanishing of the tension field $\tau (f)$ of $f$ (see for instance \cite{nishi}).

A more general context to study harmonic maps is the following. Assume we are given an isometric action of the fundamental group $\pi_{1}(M)$ of $M$ on a riemannian manifold $\widetilde{N}$ (which plays the role of the universal cover of $N$ in the previous discussion):
$$\rho : \pi_{1}(M)\to {\rm Isom}(\widetilde{N}),$$

\noindent where ${\rm Isom}(\widetilde{N})$ is the isometry group of $\widetilde{N}$. One can consider maps $f : \widetilde{M} \to \widetilde{N}$ which are equivariant with respect to the natural action of $\pi_{1}(M)$ on $\widetilde{M}$ and the action $\rho$ on $\widetilde{N}$ (here $\widetilde{M}$ is the universal cover of $M$). Such maps exist when $\widetilde{N}$ is nonpositively curved and simply connected. The energy density $$e(f): \widetilde{M}\to \mathbb{R}$$ of such an $f$ is invariant under the action of $\pi_{1}(M)$, and can be thought of as a function on $M$. Hence one can still define the energy of $f$, and once again $f$ is said to be harmonic if it is a critical point of the energy functional. If $\widetilde{N}=G/K$ is a symmetric space of noncompact type, it is known  that a harmonic equivariant map exists when the Zariski closure of the image of $\rho$ is a reductive subgroup of $G$. This is due to Corlette, see~\cite{abckt,period}.

When $M$ is K\"{a}hler, one can say a lot more about harmonic maps (this is the content of the so-called {\it nonabelian Hodge theory}).  In this context the equation of harmonic maps takes the following form. Let $\omega$ be the K\"{a}hler form on $M$. The Hodge star operator applied to the tension field gives rise to a $2n$-form $\ast \tau(f)$ on $M$ ($n={\rm dim}_{\mathbb{C}}(M)$) with values in the bundle $f^{\ast}{\rm T}\widetilde{N}$ which has the following expression:
$$\ast \tau (f)= \frac{-1}{(n-1)!}\omega^{n-1} \wedge d_{\nabla}(df\circ i),$$

\noindent where $d_{\nabla}$ is the natural operator on $f^{\ast}{\rm T}\widetilde{N}$-valued differential forms induced by the pull-back of the Levi-Civita connexion of $\widetilde{N}$ (see \cite{period}, section 14.2 page 363 for a proof).  Hence the harmonicity of $f$ is equivalent to the equation: 
$$\omega^{n-1}\wedge d_{\nabla}(df\circ i)=0.$$

We now denote again by $X$ a compact K\"{a}hler manifold, with fundamental group $\Gamma$. We fix a Zariski dense representation $\rho : \Gamma \to {\rm Isom}(\HH)$. According to the previous discussion, there exists a map $f : \widetilde{X}\to \HH$ which is $\Gamma$-equivariant and harmonic. Moreover, $f$ satisfies the following properties:
\begin{enumerate}
\item It is pluriharmonic, i.e. satisfies $d_{\nabla}(df\circ i)=0$.
\item The $(0,1)$ part of $d_{\nabla}$, denoted by $d_{\nabla}^{0,1}$, satisfies $$\left( d_{\nabla}^{0,1}\right)^{2}=0.$$ According to the Koszul-Malgrange integrability theorem \cite{km}, this implies that there exists a holomorphic structure on the bundle $f^{\ast}{\rm T}\HH \otimes \mathbb{C}$ for which  $\overline{\partial}=d_{\nabla}^{0,1}$. 
\item The complex linear part $\alpha := df^{1,0}$ of the differential of $f$ is a holomorphic $1$-form with values in $f^{\ast}{\rm T}\HH \otimes \mathbb{C}$ (for the previous holomorphic structure). 
\item The complex rank of $\alpha$ is everywhere $\le 1$, consequently the real rank of $f$ is everywhere $\le 2$. 
\end{enumerate}

The first three results are proved using the classical Bochner-Siu-Sampson formula (see \cite{abckt} or \cite{period}). The last one is due to Sampson \cite{sampson2}. Note that the first three results hold more generally when the real hyperbolic space is replaced by any symmetric space of noncompact type. The last one on the other hand is very specific to $\HH$. 

All the notions that we have just described still make sense when the representation $\rho : \Gamma \to \SO$ is replaced by a nonelementary action of $\Gamma$ on the space $\HHI$ and $f$ is an equivariant map from $\widetilde{X}$ to $\HHI$. The existence of a harmonic map in this case is a result of Korevaar and Schoen~\cite{korschoen0,korschoen}. In the next two subsections, we will describe this result and explain why the four properties listed above still hold in this context. 

\subsection{Korevaar and Schoen's harmonic map}\label{resultat-ks}

 We start with a remark. We need to do some differential geometric calculations in $\HHI$. To this end, observe that the Levi-Civita connection of $\HHI$ can be defined in the usual way as the unique torsion-free connection on the tangent bundle ${\rm T}\HHI$ of $\HHI$, which is compatible with the metric. When we use the Poincar\'e ball model for $\HHI$, i.e. we identify $\HHI$ with the unit ball $\mathscr{B}\subset \mathscr{H}$, endowed with the metric $\frac{4\langle \cdot , \cdot \rangle}{(1-\vert x\vert^{2})^{2}}$ we see easily that if $X : \mathscr{B}\to \mathscr{H}$ is a smooth map thought of as a vector field on $\HHI$, one has: 
$$\nabla X=dX+ A(X),$$
\noindent where $A$ is a smooth $1$-form with values in the space of endomorphisms of $\mathscr{H}$. One can compute easily the precise expression of $A$ in terms of the function $\frac{4}{(1-\vert x\vert^{2})^{2}}$ and of its gradient, but we will not need it. In the following if $v$ is a vector in ${\rm T}\HHI$, we will write $A_{v}$ for the value of the $1$-form $A$ on $v$, so that $A_{v} : \mathscr{H}\to \mathscr{H}$ is a continuous linear map.

We now fix a nonelementary homomorphism $\rho : \Gamma \to {\rm Isom}(\HHI)$ and we want to study $\rho$ via an associated harmonic map. The existence of an equivariant harmonic map in this context is a consequence of the work of Korevaar and Schoen~\cite{korschoen0,korschoen}, see Theorem~2.3.1 in~\cite{korschoen}. The minimizing map constructed by Korevaar and Schoen is Lipschitz and belongs to a certain Sobolev space defined in metric terms, see \cite{korschoen0} for the definition. Since in our case the target is a smooth infinite-dimensional manifold, we will see that the harmonic map is in fact smooth. It is well-known that harmonic maps are smooth once we know their continuity (see for instance \cite{borchers} for a proof for maps from the Euclidean space into the sphere, as well as \cite{schoenuhl} for a discussion of the regularity of harmonic maps). In our case the stronger fact that $f$ is Lipschitz was established in \cite{korschoen}; deducing the smoothness from this is easy. Since we were not able to find a reference for this in the litterature, we now prove this fact. 

In the following $B$ will be a closed ball in $\widetilde{X}$ contained in some coordinate chart, $g$ will denote the riemannian metric on $\widetilde{X}$ and $\Delta_{g}$ the Laplacian associated to $g$. We will need to deal with various Lebesgue and Sobolev spaces of maps with values in the real separable Hilbert space $\mathscr{H}$. One defines the space $L^{p}(B,\mathscr{H})$ (for $p\ge 1$), that we will simply denote by $L^{p}$, as the space of all measurable maps $u : B \to \mathscr{H}$ such that $\int_{B}\vert \vert u(x)\vert \vert^{p}dx<+\infty$ (where $dx$ stands for the riemannian measure associated to $g$ and $\vert \vert \cdot \vert \vert$ is the norm on $\mathscr{H}$). We define the Sobolev spaces $W^{k,p}(B,\mathscr{H})$ (or $W^{k,p}$ for short) in the usual way, using the scalar product on $\mathscr{H}$. Finally, we will write $u\in L^{p}_{loc}$ (or $u\in W^{1,p}_{loc}$ or $u\in W^{2,p}_{loc}$...) if the restriction of $u$ to any closed ball contained in the interior of $B$ is in $L^{p}$. See \cite{diestel} for a detailed study of Lebesgue spaces of vector-valued functions. Note that in the present case, since we are dealing with functions with values in a Hilbert space, the integral of a measurable map $u : B \to \mathscr{H}$ such that $\int_{B}\vert \vert u\vert \vert dx<+\infty$ can be defined easily by duality.

We start with an easy proposition. 

\begin{prop}\label{lemme-reg} Assume that the measurable map $u : B \to \mathscr{H}$ is in $W^{k,2}$ for any integer $k$. Then $u$ is of class $C^{\infty}$.
\end{prop}
{\it Proof.} Writing $(u_{i})_{i\ge 0}$ for the coordinates of $u$ in a Hilbert basis of $\mathscr{H}$, the hypothesis implies that each $u_{i}$ is in the usual Sobolev space $W^{k,2}$ and for each fixed $k$, the series $\sum_{i}\vert u_{i}\vert_{W^{k,2}}^{2}$ is convergent. Using the classical regularity theorems each $u_{i}$ is smooth and for each $r$, and each $k$ large enough compared to $r$, there exists a constant $D_{r,k}$ such that $\vert u_{i}\vert_{C^{r}}\le D_{r,k}\vert u_{i}\vert_{W^{k,2}}$. Thus, the series $\sum_{i}\vert u_{i}\vert_{C^{r}}^{2}$ is convergent, which implies that $u$ is $C^{\infty}$. 
\hfill $\Box$

 The following proposition is classical, at least for real-valued functions. We will briefly remind its proof in section~\ref{analyse}, to explain why it extends to the case of vector-valued maps.

\begin{prop}\label{laplacelp} If a map $u$ is in $W^{1,p}_{loc}$ for all $p$ and if $\Delta_{g}u$ is in $L^{p}_{loc}$ for all $p$, then $u$ is in $W^{2,p}_{loc}$ for all $p$. 
\end{prop}

 We are now ready to prove the:

\begin{prop}The harmonic map $f : \widetilde{X} \to \HHI$ is of class $C^{\infty}$. 
\end{prop}
{\it Proof.} We fix a point $x_{0}\in \widetilde{X}$ and work in a ball $B$ around  $x_{0}$. We denote by $x_{1}, \ldots , x_{2n}$ the coordinates on $B$. We also use the Poincar\'e ball model for $\HHI$ i.e. we identify $\HHI$ with the unit ball $\mathscr{B}$ of $\mathscr{H}$ endowed with the metric $\langle \cdot , \cdot \rangle_{hyp}:= \frac{4\langle \cdot , \cdot \rangle}{(1-\vert x\vert^{2})^{2}}$. Since $f$ is continuous, we can assume, by shrinking $B$ if necessary, that $f(B)$ is contained in a small ball of $\HHI$ where the hyperbolic and the flat metric are comparable. This implies that $f$ (restricted to $B$) is in the Sobolev space defined by Korevaar and Schoen for any of the two metrics. But for the flat metric on $\mathscr{H}$, Korevaar and Schoen's Sobolev space coincides with the usual space $W^{1,2}$ defined above. This is Theorem 1.6.2 in \cite{korschoen0} (that theorem deals with real-valued functions but can be adapted easily to the case where the target is a Hilbert space).

Observe that since $f$ is Lipschitz, the partial derivatives of order $1$ of $f$ are in fact $L^{\infty}$. As in \cite{borchers}, we are going to show that for any integer $k$ and any multiindex $I=\{\alpha_{1},\ldots , \alpha_{k}\}$ of length $k$, the partial derivative $\partial_{I}f:=\partial_{\alpha_{1}}\cdots \partial_{\alpha_{k}}f$ is in $L^{p}_{loc}$ for any $p>1$. Together with Proposition~\ref{lemme-reg}, this will prove the proposition.

First observe that the fact that $f$ minimizes energy implies (as in the smooth case) that for any smooth compactly supported map $V : B \to \mathscr{H}$ (thought of as a vector field along $f$) one has:
$$\int \langle \nabla V , df\rangle_{hyp} dx=0.$$

\noindent Here the symbol $dx$ represents the integration on $B$ with respect to the volume form induced by $g$ and the scalar product of the two $1$-forms is defined in the usual way. When $f$ is smooth this last integral equals $-\int \langle V , \tau(f) \rangle_{hyp} dx$ where $\tau (f)$ is the tension field of $f$ (see \cite{nishi}, page 102). Remembering that $\nabla V=dV+A(V)$, it is easy to deduce from the previous equality that $f$ satisfies (in the sense of distributions) the equation:

$$\Delta_{g}f+\sum_{i,j=1}^{2n}g^{ij}A_{\frac{\partial f}{\partial x_{i}}}\left(\frac{\partial f}{\partial x_{j}}\right)=0\;\;\;\;\;\; (\ast)$$

\noindent We now prove by induction on $k$ that $\partial_{I}f\in L^{p}$ for any $p$ and for any multiindex $I$ of length $k$. This is true for $k=1$. Assuming it is true for a certain integer $k\ge 1$ we prove it for $k+1$. Let $I$ be a multi-index of length $k-1$. According to $(\ast)$, the function $\Delta_{g}\partial_{I}f$ is a linear combination of products of partial derivatives of length at most $k$ of $A_{i}(f)$, $g^{ij}$ and $f$ and of the term $\Delta_{g}\partial_{I}f-\partial_{I}\Delta_{g}f$. This last term involves derivatives of $f$ up to the order $k$. By induction all these terms are in $L^{p}$ for any $p$, hence by the H\"{o}lder inequality $\Delta_{g}\partial_{I}f$ is in $L^{p}$ for any $p$. By Proposition~\ref{laplacelp}, $\partial_{I}f\in W^{2,p}$, hence $f\in W^{k+1,p}$.\hfill $\Box$

Once $f$ is known to be smooth, we can define as before the operator $d_{\nabla}$ acting on differential forms with values in the bundle $f^{\ast}{\rm T}\HHI \to \widetilde{X}$, and one checks as in the finite-dimensional case that the harmonicity of $f$ is equivalent to the usual equation:
$$\omega^{n-1}\wedge d_{\nabla}(df\circ i)=0.$$

\subsection{The Bochner-Siu-Sampson formula}

The next step is to check that the usual Bochner-Siu-Sampson formula still holds in this context, showing that $f$ is pluriharmonic and that the curvature $\left( d_{\nabla}\right)^{2}$ is purely of type~$(1,1)$.

\begin{prop} The following holds:
\begin{enumerate}
\item The map $f$ is pluriharmonic i.e. $d_{\nabla}(df\circ i)=0$,
\item The pull-back of the curvature of $\HHI$ (thought of as a $2$-form with values in ${\rm End}(T\HHI)$ is of type $(1,1)$ which implies that $D:=d_{\nabla}^{1,0}$ satisfies $D^{2}=0$. 
\item $D(d^{1,0}f)=0$.
\end{enumerate}
\end{prop}
{\it Proof.} The proof of this proposition in the finite-dimensional case is an application of the classical Bochner-Siu-Sampson formula, see~\cite{abckt} or \cite{period}, section 14.2. Here we only remind the reader of the strategy of the proof and explain why it goes through without changes when the target is an infinite-dimensional hyperbolic space. 

We first introduce the following notation. If $\beta$ and $\delta$ are differential forms on $\widetilde{X}$ with values in the bundle $f^{\ast}{\rm T}\HHI \otimes \C$ one defines naturally the $k+l$ form $\langle \beta\wedge \delta \rangle$:  if $\beta=\sum \beta_{i}s_{i}$ and $\delta=\sum_{j}\delta_{j}v_{j}$ (where the $\beta_{i}$ and the $\delta_{j}$ are complex valued forms and the $s_{i}$ and the $t_{j}$ are sections of $f^{\ast}{\rm T}\HHI\otimes \C$) we define:
$$\langle \beta \wedge \delta\rangle:=\sum_{i,j}\beta_{i}\wedge \overline{\delta_{j}}\langle s_{i},t_{j}\rangle$$
\noindent where $\langle s_{i},t_{j}\rangle$ is computed using the natural hermitian metric on $f^{\ast}{\rm T}\HHI \otimes \C$. 

The proof consists in showing that the real-valued $2n$-form (which is well-defined on $X$)

$$\langle d_{\nabla}^{2}(df\circ i)\wedge df\circ i \rangle \wedge \omega^{n-2}+\langle d_{\nabla}(df\circ i)\wedge d_{\nabla}(df\circ i) \rangle \wedge \omega^{n-2}$$

\noindent is at the same time exact and pointwise nonpositive. Stokes theorem then implies that it is identically zero. The exactness of this form is proved as in the finite-dimensional case. As for the sign, the harmonic equation tells us that for any section $s$ of $f^{\ast}{\rm T}\HHI$, the $2$-form $\langle d_{\nabla}(df\circ i),s\rangle \wedge \omega^{n-1}$ is primitive. This implies that the second term in the formula above is nonpositive and vanishes only if $d_{\nabla}(df\circ i)=0$ identically. The first term is also nonpositive and vanishes precisely if the curvature $d_{\nabla}^{2}$ is of type $(1,1)$. The proof of this fact goes through exactly as in the finite-dimensional case, see \cite{abckt}, page 74. 

The proof of the third point of the proposition is now easy. The $1$-form $$\alpha:=df^{1,0}=df-idf\circ i$$ is $d_{\nabla}$-closed since both $df$ and $df\circ i$ are. Since $\alpha$ is of type $(1,0)$ this implies that both $d_{\nabla}^{1,0}\alpha$ and $d_{\nabla}^{0,1}\alpha$ are $0$.\hfill $\Box$

In the following we will continue to denote by $\alpha$ the complex linear part of the differential of $f$: $$\alpha:=df^{1,0}.$$ \noindent Exactly as in~\cite{sampson2}, the previous proposition implies that the complex rank of $\alpha$ is at most $1$. We will briefly summarize the proof of this fact here. See~\cite{abckt}, Chapter 6, page 77, for more details. We denote by $R$ the curvature tensor of $\HHI$ and by $R_{\mathbb{C}}$ its complex linear extension to $\left( {\rm T}\HHI \otimes \mathbb{C}\right)^{\otimes 4}$. Denoting by $\langle \cdot , \cdot \rangle$ the complex bilinear extension of the Riemannian metric of $\HHI$ to ${\rm T}\HHI \otimes \C$, one has:
$$R_{\mathbb{C}}(X,Y,Z,W)=\langle Y,Z\rangle \langle X , W \rangle-\langle X , Z \rangle \langle Y , W \rangle.\;\;\;\;\;(\ast)$$
\noindent Now, one can show (see \cite{abckt}) that the second item of the previous proposition, namely, the fact that the pull-back by $f$ of the curvature is of type $(1,1)$ is equivalent to the fact that the {\it Hermitian sectional curvature} is zero on the image of $\alpha$: for all $X$, $Y$ in the image of $\alpha$ we have:
$$R_{\C}(X,Y,\overline{X},\overline{Y})=0.$$
\noindent From the expression $(\ast)$ for $R_{\C}$, one sees easily that this implies that the image of $\alpha$ has complex dimension at most $1$. This implies that the real rank of $f$ is not greater than $2$.

The next step is to check that the $1$-form $\alpha$ is actually holomorphic for a suitable holomorphic structure on the bundle $f^{\ast}{\rm T}\HHI\otimes \C$. To achieve this, we will need an analogue of the Koszul-Malgrange theorem in this infinite dimensional context. The statement of this theorem for a Hilbert space bundle does not present any difficulty. Knowing that the operator $D$ acting on sections of the bundle $f^{\ast}{\rm T}\HHI \otimes \C \to \widetilde{X}$ satisfies $D^{2}=0$, we want to prove that there exist local trivializations of this bundle in which the operator $D$ becomes the usual $\overline{\partial}$ operator. We will see in section~\ref{section-dbar} that the proof in the finite dimensional case can be adapted without difficulty to this context. Hence we assume from now on that such trivializations exist. Observe however that around any point $x$ such that $\alpha_{x}\neq 0$, the image of $\alpha$ is a rank $1$ subbundle of $f^{\ast}{\rm T}\HHI$ which is stable by the operator $D$. Therefore, we can apply the classical Koszul-Malgrange theorem to this rank $1$ bundle around any point where $\alpha$ is not zero. So we only need the infinite-dimensional version of the theorem to study points where $\alpha$ vanishes.

%%%%%%%%%%%%%%%%%%%%%%%%%%%%%%%%%%%%%%%%%%%%%%%%%%%%%%%%%%%%%%%%%%%%%%%%%%%%%%%%%%%%%%%
%%%%%%%%%%%%%%%%%%%%%%%%%%%%%%%%%%%%%%%%%%%%%%%%%%%%%%%%%%%%%%%%%%%%%%%%%%%%%%%%%%%%%%%

\section{Actions on infinite dimensional real hyperbolic spaces}\label{section-infinite}

In this section, we prove Theorem~\ref{infinite}. We consider a minimal non-elementary action $\Gamma \to {\rm Isom}(\HHI)$ and work with the associated equivariant harmonic map $f : \widetilde{X} \to \HHI$. First, let us observe that one can assume that the rank of $f$ equals $2$ on a dense open set. Otherwise the image of $f$ would be a geodesic or a point (see \cite{sampson}, Theorem 3, the proof works just as well in our context) and the representation $\rho$ would be elementary.

\subsection{Factorization}\label{factorization}

The goal of this subsection is to prove the following theorem. It is classical in spirit and similar to many results in the litterature. Although a very similar statement appears in \cite{jostyau}, we include a complete proof here.

\begin{main}\label{prop-factor-infinite}There exists a holomorphic map $\pi : \widetilde{X}\to \Delta$ and a harmonic map \mbox{$u : \Delta \to \HHI$} such that:
$$f = u\circ \pi.$$
\noindent Moreover there exists a homomorphism $\theta : \Gamma \to \PSLR$ and a homomorphism $$\Psi : \Lambda:=\theta(\Gamma)\to {\rm Isom}(\HHI)$$ \noindent such that $\rho=\Psi \circ \theta$. The map $\pi$ is $\theta$-equivariant and the map $u$ is $\Psi$-equivariant. 
\end{main}

To prove this theorem, we follow an idea from \cite{mok88}, section 2.2. In the course of the proof, we will use the fact that the map $f$ is real-analytic (as a solution of an elliptic equation with real-analytic coefficients), as well as the following lemma, whose proof we omit. 

\begin{lemma} Let $O$ be a manifold and let $R\subset O\times O$ be the graph of an equivalence relation on $O$. If $R$ is closed and if the two projections $p_{1},p_{2} : R\to O$ are open maps ($R$ being endowed with the topology induced from $O\times O$), then the quotient space $O_{R}$ is Hausdorff and the map  $O\to O_{R}$ is an open map.
\end{lemma}

Let $U$ be the open set of $\widetilde{X}$ where the holomorphic $1$-form $\alpha$ is nonzero. Note that its kernel defines a holomorphic foliation $\mathscr{F}$ on $U$ (for instance because the kernel of $\alpha$ coincides with that of $df$ on a dense open set). The foliation $\mathscr{F}$ extends to $O:=\widetilde{X}-Z$ where $Z$ is an analytic set of codimension at least $2$ in $\widetilde{X}$, called the singular set of $\mathscr{F}$. To prove Theorem~\ref{prop-factor-infinite}, it is enough to construct a holomorphic equivariant map $O\to \Delta$, since such a map will automatically extend to $\widetilde{X}$.

We define a subset $R$ of $O\times O$ as follows. The pair $(x,y)\in O\times O$ is in $R$ if we can find two germs of holomorphic maps $\varphi_{x} : (\C,0)\to (O,x)$ and $\varphi_{y} : (\C,0)\to (O,y)$ whose images are not contained in a leaf of $\mathscr{F}$ and such that $f\circ \varphi_{x}=f\circ \varphi_{y}$. Let us make some comments about this definition. 

\begin{enumerate}
\item First, if one takes coordinates $z=(z_{1}, z')\in \Delta(0,\varepsilon)\times \Delta(0,\delta)\subset \C \times \C^{n-1}$ and $w=(w_{1},w')\in \Delta(0,\varepsilon)\times \Delta(0,\delta)\subset \C \times \C^{n-1}$ centered at $x$ and $y$ respectively, and such that $\mathscr{F}$ is locally defined by the equations $dz_{1}=0$ and $dw_{1}=0$, one sees that $R$ is locally the ``pull-back" by the map $(z,w)\mapsto (z_{1},w_{1})$ of a subset $R_{0}$ of $\Delta(0,\varepsilon) \times \Delta(0,\varepsilon)$. Hence, by projecting onto the $z_{1}$ (or $w_{1}$) coordinate we can think of the two maps $\varphi_{x}$ and $\varphi_{y}$ as maps $(\C,0)\to (\C,0)$. We will always do this in what follows, and work with $R_{0}$. Observe also that in such coordinates one has $f(z_{1},z')=f(z_{1})$ (and similarly $f(w_{1},w')=f(w_{1})$).

\item Second, we observe that $R$ defines an equivalence relation on $O$. The only nontrivial point is the fact that $R$ defines a transitive relation but this essentially follows from the fact that two (germs of) ramified covers of $(\C,0)$ share a common cover. 

\item Third, if $(x,y)\in R$, Let $t >0$ be such that $\varphi_{x}$ and $\varphi_{y}$ are both defined on the disc $\Delta(0,t)$. Then the germ of complex analytic set 
$$\{(\varphi_{x}(z),\varphi_{y}(z)), z\in \Delta(0,t)\}\subset R_{0}\subset \C^{2}$$

\noindent is contained in the real analytic set $R_{1}:=\{(z,w),f(z)=f(w)\}\subset \C^{2}$. Since $f$ is not constant on any open set, $R_{1}$ cannot contain a germ of analytic curve whose projection to one of the two factors of $\C^{2}$ is constant. Hence, we see that $R_{0}$ is exactly the set of points of $R_{1}$ which are contained in a positive dimensional germ $A$ of complex analytic set with $A\subset R_{1}$. According to \cite{dm} this set is closed. Hence the graph $R$ of the equivalence relation we just defined is closed. The projection maps $p_{1},p_{2} : R\subset O\times O\to O$ are open. According to the lemma mentioned above the quotient $O_{R}$ of $O$ by the equivalence relation $R$ is Hausdorff and the projection map $\pi : O \to O_{R}$ is open.  
\end{enumerate}

We now prove that $O_{R}$  is a Riemann surface. We start with the following:

\begin{prop}\label{propcomp} For any point $p$ in $O$ there is a neighborhood $U$ of $(p,p)$ in $O\times O$ such that $R\cap U$ is a complex analytic set. 
\end{prop}

\noindent {\it Proof.} We choose coordinates $(z_{1},\ldots , z_{n})=(z_{1},z')$ in a polydisc $$V=\Delta(0,\varepsilon)\times \Delta (0,\delta)\subset \C\times \C^{n-1}$$ centered at $p$ in such a way that the foliation $\mathscr{F}$ is defined on $V$ by the equation $dz_{1}=0$. Let $D$ be the disc $\Delta(0,\varepsilon)\times \{0\}$. In these coordinates, the harmonic map depends on $z_{1}$ only. After choosing a trivialization of the bundle $f^{\ast}{\rm T}\HHI\otimes \C$ near $p$, the form $\alpha$ can be thought of as a $1$-form depending on $z_{1}$ only, with values in a fixed complex Hilbert space ${\bf H}$. Let $(x,y)\in V \times V \cap R$ and let $\varphi_{x},\varphi_{y} : (\C,0)\to D$ be the associated maps as in the definition of $R$. From the equation $f\circ \varphi_{x}=f\circ \varphi_{y}$ one gets $\varphi_{x}^{\ast}\alpha=\varphi_{y}^{\ast}\alpha$. By choosing a basis $(e_{1}, e_{2},\ldots)$ of ${\bf H}$ and maybe making a change of coordinate on $z_{1}$, we can assume that:
$$\alpha= z_{1}^{k}dz_{1}e_{1}+\sum_{j\ge 2}a_{j}(z_{1})dz_{1}e_{j},$$

\noindent for some integer $k$. By projecting the equation $\varphi_{x}^{\ast}\alpha=\varphi_{y}^{\ast}\alpha$ on the first coordinate, we get:
$$\varphi_{x}^{k}\varphi_{x}'=\varphi_{y}^{k}\varphi_{y}',$$

\noindent Hence the set $R_{0}\subset D\times D$ (which was defined to be the image of $R$ by the map $(z,w)\mapsto (z_{1},w_{1})$) is contained in the set
$$\bigcup_{c\in \C}A_{c},$$

\noindent where $A_{c}:=\{(z,w)\in D\times D, z^{k+1}-w^{k+1}=c\}$. We are going to see that, up to replacing $D$ by a smaller disc, one has $R_{0}\subset A_{0}$. This will prove the proposition. If this is false, there is a sequence $c_{n}$ of nonzero complex numbers converging to $0$ such that $$A_{c_{n}} \cap R_{0}\neq \emptyset . \;\;\;\;\;\; (\ast)$$

\noindent We will prove that this implies that the rank of $f$ equals $1$ on an open set, a contradiction. For $n$ large enough the set $A_{c_{n}}$ is a connected Riemann surface in $D\times D$. If it intersects $R_{0}$, the two harmonic maps $f\circ p_{1}$ and $f\circ p_{2}$ (where the $p_{i}$ are the projections $D\times D\to D$) coincide on an open set of $A_{c_{n}}$. Hence they are equal on $A_{c_{n}}$ and $A_{c_{n}}$ is entirely contained in $R_{0}$. We thus assume that 
$$A_{c_{n}}\subset R_{0}\subset R_{1},$$
\noindent for all $n$. Now, for any point $x\in D$ one can find a sequence of points $x_{n}$ such that 
$$x^{k+1}-x_{n}^{k+1}=c_{n}$$
\noindent and $x_{n}\to x$. Since $A_{c_{n}}\subset R_{1}:=\{(a,b), f(a)=f(b)\}$ we have $f(x)=f(x_{n})$, and by taking a limit point of the sequence $\frac{x-x_{n}}{\vert x-x_{n}\vert}$ one sees that ${\rm rank}\, df_{x}\le 1$ for $x\in D$. This is the contradiction we were looking for. \hfill $\Box$

With the help  of the previous proposition, we can conclude the proof. The quotient of a polydisc $V$ as in the proposition above naturally has the structure of a Riemann surface. Since $\pi$ is an open map, this endows the open set $\pi(V)\subset O_{R}$ with a complex structure. It is then easy to check that the various complex structures on the open sets $\pi(V)$ (for various choices of $V$) are compatible and make the map $\pi$ holomorphic.   

\begin{rem} The idea of choosing a complex analytic equivalence relation whose graph is contained in the real analytic set $\{(a,b),f(a)=f(b)\}$ is due to Mok, see \cite{mok88}. However, we use it in a slighty different way. Note also that we cannot prove immediately that $R$ is complex analytic everywhere, this is why in Proposition~\ref{propcomp} we restrict to a neighborhood of the diagonal. Indeed, if we start with a point $(p,q)\in R$ ($p\neq q$), there is no way to identify ({\it a priori}) the fibers of the bundle $f^{\ast}{\rm T}\HHI \otimes \C$ near $p$ and near $q$. 
\end{rem}

From now on we will denote by $\Sigma$ the Riemann surface $O_{R}$. At this stage we have a holomorphic map $\pi : \widetilde{X}\to \Sigma$ and a continuous map $u : \Sigma \to \HHI$ such that $f=u\circ \pi$. We still have to justify that $u$ is harmonic. But this is clear in a neighborhood of any point which is not a critical value of $\pi$ since $f$ is pluriharmonic. To establish the smoothness and the harmonicity of $u$ near a critical value of $\pi$, we borrow an argument from \cite{mok}, page 577. Let $p$ be a critical value of $\pi$ and let $x\in \widetilde{X}$ be such that $\pi (x)=p$. We can assume that $x$ is a regular point of the foliation defined by the kernel of $\alpha$. We choose coordinates $z=(z_{1}, \ldots , z_{n})$ centered at $x$ such that $\pi$ is locally represented by the map $z\mapsto z_{1}^{k}$ for some integer $k$. From the equation 
$$f(z)=u(z_{1}^{k}),$$
\noindent we want to deduce that $u$ is smooth and harmonic. We will denote by $D_{0}$ a small disc centered at the origin in the complex plane (corresponding to the $z_{1}$ coordinate) and by $D$ its image under the map $z_{1}\mapsto z_{1}^{k}$. We also write $f(z)=f(z_{1})$ in what follows. Let $v : D \to \HHI$ be the harmonic map whose boundary value is equal to $u$. It exists according to~\cite{korschoen0}, Theorem 2.2.  Now the map $\widetilde{v}$ defined on $D_{0}$ by $\widetilde{v}(z_{1})=v(z_{1}^{k})$ is harmonic and equals $f$ on the boundary of $D_{0}$. Hence it must equals $f$ on $D_{0}$. This implies that $u=v$ on $D$ and proves the result. 

We now explain why we can assume that $\Sigma$ is the unit disc $\Delta$. First, by the maximum principle, the Riemann surface $\Sigma$ cannot be compact. Then, we observe that the map $\pi : \widetilde{X}\to \Sigma$ is equivariant with respect to a homomorphism $\theta : \Gamma \to {\rm Aut}(\Sigma)$. The group ${\rm Aut}(\Sigma)$ cannot be solvable as a consequence of the following lemma.

\begin{lemma} Let $G$ be a finitely generated solvable subgroup of ${\rm Isom}(\HHI)$. Then either $G$ stabilizes a geodesic of $\HHI$ or $G$ fixes a point in $\HHI \cup \partial \HHI$.
\end{lemma}

 Note that this lemma is very general and remains true for groups of isometries of $\delta$-hyperbolic spaces. A proof in the case of $\HHI$ can be found in \cite{cantat}, section 6.3, for instance. It implies that the universal cover of $\Sigma$ is isomorphic to the unit disc. Although $\Sigma$ might not be simply connected at the beginning, one can always lift the map $\pi$ (as well as the representation $\Gamma \to {\rm Aut}(\Sigma)$) to the unit disc. To be also able to lift the map $u : \Sigma \to \HHI$ to a map defined on $\Delta$, it is enough to check that the lift $\widetilde{\pi} : \widetilde{X}\to \Delta$ of $\pi$ is surjective. Since the image of $\widetilde{\pi}$ is open we need to check that it is also closed. If the automorphism group of $\Sigma$ is discrete, it is easy to see that the image of $\widetilde{\pi}$ is closed. We conclude by observing that the only hyperbolic Riemann surface whose automorphism group is not discrete and not solvable is the unit disc itself; if $\Sigma = \Delta$, $\widetilde{\pi}=\pi$ and we know from the beginning that this map is surjective. Hence, up to replacing $\Sigma$ by its universal cover and $\pi$ by a lift, we now assume that $\Sigma=\Delta$.    
  
Before concluding the proof of Theorem~\ref{prop-factor-infinite}, we state the following lemma:

\begin{lemma}\label{unique} If an isometry $g : \HHI \to \HHI$ fixes pointwise the image of the map $f$ then $g$ is the identity.
\end{lemma}
{\it Proof.} Let $W$ be the intersection of all closed totally geodesic subspaces of $\HHI$ containing the image of $f$. It is itself a closed totally geodesic subspace of $\HHI$, which is invariant by the action of $\Gamma$. Since we assumed that the action was minimal, we have $W=\HHI$. We conclude the proof by observing that the fixed point set of $g$ contains $W$. \hfill $\Box$

We can now conclude. Let $\Lambda$ be the image of $\theta$. According to Lemma~\ref{unique}, the kernel of $\theta$ is contained in the kernel of $\rho$. So we can write $\rho= \Psi \circ \theta$, where $\Psi : \Lambda \to {\rm Isom}(\HHI)$ is a homomorphism. This finishes the proof of Theorem~\ref{prop-factor-infinite}. 
 
\subsection{Conclusion of the proof of Theorem~\ref{infinite}}\label{conclusionoftheproof}
 If the group $\Lambda$ introduced above is discrete in $\PSLR$, the representation factors through a fibration onto a hyperbolic $2$-orbifold. Indeed one gets a map $X\to \Delta/\Lambda$, which might not have connected fibers. However, we can apply the Stein factorization theorem to it to obtain a fibration of $X$ onto a hyperbolic $2$-orbifold. The obvious, but crucial observation, which appears already in \cite{gs} (Lemma 9.4), is now the following: {\it if $\Lambda$ is not discrete, it is topologically dense in $\PSLR$}.  Assuming that $\Lambda$ is not discrete, we now prove that the homomorphism $\Psi : \Lambda \to {\rm Isom}(\HHI)$ extends continuously to $\PSLR$. Let $(\gamma_{n})$ be a sequence in $\Lambda$ which converges to $g\in \PSLR$. We are going to prove that $\rho (\gamma_{n})$ converges to an isometry of $\HHI$. First observe the following:

\begin{lemma} Let $(g_{n})$ be a sequence of isometries of $\HHI$. Then the set 
$$\{x\in \HHI, (g_{n}(x)) \; is \; convergent\}$$
is a closed totally geodesic space. 
 \end{lemma}
 {\it Proof.} The fact that this set is closed follows from the completeness of $\HHI$ and the equicontinuity of the sequence $(g_{n})$. The fact that it is totally geodesic follows from the fact that the geodesic ray passing through two given points of $\HHI$ depends continuously on these two points. \hfill $\Box$

The two sequences $\rho(\gamma_{n})$ and $\rho(\gamma_{n})^{-1}$ are pointwise convergent on the image of the harmonic map. According to the lemma and to the fact that the smallest totally geodesic subspace of $\HHI$ containing ${\rm Image}(f)$ is $\HHI$, they are pointwise convergent on $\HHI$. This proves that the homomorphism $\Psi$ extends to a continuous homomorphism $\PSLR \to {\rm Isom}(\HHI)$. This concludes the proof of Theorem~\ref{infinite}: if the representation $\rho : \Gamma \to {\rm Isom}(\HHI)$ does not factor through a fibration onto a hyperbolic $2$-orbifold, it factors through a continuous homomorphism from $\PSLR$ to ${\rm Isom}(\HHI)$. 

As mentioned in remark~\ref{symspace}, the same arguments prove that a Zariski dense representation of $\Gamma$ in a noncompact simple Lie group with trivial center and different from $\PSLR$ factors through a hyperbolic $2$-orbifold as soon as the associated harmonic map $$f : \widetilde{X}\to G/K$$ \noindent satisfies ${\rm rank}_{\mathbb{C}}\, df^{1,0}=1$. Indeed this hypothesis ensures that there is a holomorphic foliation of codimension $1$ on $\widetilde{X}$ and we can argue as before. We use the fact that there is no nontrivial subsymmetric space containing the image of $f$; as a consequence any sequence in $G$ which converges pointwise on the image of $f$ has to converge everywhere. 

%%%%%%%%%%%%%%%%%%%%%%%%%%%%%%%%%%%%%%%%%%%%%%%%%%%%%%%%%%%%%%%%%%%%%%%%%%%%%%%%%%%%%%%
%%%%%%%%%%%%%%%%%%%%%%%%%%%%%%%%%%%%%%%%%%%%%%%%%%%%%%%%%%%%%%%%%%%%%%%%%%%%%%%%%%%%%%%

\section{K\"{a}hler groups in the Cremona group}\label{section-cremo}

Here, we recall the construction of the {\it Picard-Manin space} on which the Cremona group acts faithfully, before proving Theorem~\ref{cremo}. For more details on the notions presented here, the reader can consult \cite{bfj,cantat,cantatlamy,favre,manin}.  

\subsection{The Picard-Manin space}

We consider the set $\mathscr{E}$ of all pairs $(X,p)$ where $X$ is a smooth rational surface and \mbox{$p : X\to \mathbb{P}^{2}$} is a birational morphism. Observe that if $(X,p)\in \mathscr{E}$, the map $p$ is a composition of a certain number $k$ of blow-ups and the intersection form $(\cdot , \cdot)_{X}$ on the group $H^{2}(X,\Z)$ has signature $(1,k)$ (i.e. $H^{2}(X,\Z)$ with its intersection form is isomorphic to $\Z^{k+1}$ with the form $x_{1}^{2}-x_{2}^{2}-\cdots -x_{k+1}^{2}$). Although one can define the ``action" $f^{\ast}$ of a birational map $f : X \dashrightarrow X$ on the group $H^{2}(X,\Z)$ (see \cite{dillerfavre}), the corresponding map $f\mapsto (f^{-1})^{\ast}$ is not a group homomorphism. To solve this problem, we look at the action of the Cremona group on a suitably defined limit of all cohomology groups $H^{2}(X,\Z)$ for $(X,p)\in \mathscr{E}$. The construction of this limit is due to Manin \cite{manin}.

Formally, one proceeds as follows. If $(X,p)$ and $(Y,q)$  are elements of $\mathscr{E}$, we will say that a morphism $u : Y\to X$ is admissible if $q=p\circ u$. Note that there is at most one admissible morphism between $(X,p)$ and $(Y,q)$. There is one exactly when the birational map $p^{-1}\circ q$ is holomorphic. If $u : Y\to X$ is admissible, the map $$u^{\ast} : (H^{2}(X,\Z),(\cdot , \cdot)_{X})\to (H^{2}(Y,\Z),(\cdot , \cdot)_{Y})$$ is an isometric embedding. Observe also that given $(X,p)$ and $(Y,q)$ in $\mathscr{E}$, one can always find a third surface $(Z,r)\in \mathscr{E}$ for which there exist admissible morphisms $Z\to X$ and $Z\to Y$. This allows one to define the space

$$\ZZ=\underset{\rightarrow}{{\rm lim}}\; H^{2}(X,\Z),$$

\noindent which is the inductive limit of the groups $H^{2}(X,\Z)$ (for $(X,p)\in \mathscr{E}$) with respect to all admissible morphisms. The abelian group $\ZZ$ is endowed with an intersection form $(\cdot , \cdot)_{\mathbb{P}^{2}}$. One defines an action of the group $\bir$ on $\ZZ$ as follows. Let $f : \mathbb{P}^{2} \dashrightarrow \mathbb{P}^{2}$ be a birational map. If $\alpha\in H^{2}(X,\Z)$  (for some $(X,p)\in \mathscr{E}$) and $[\alpha]$ denotes the image of $\alpha$ in $\ZZ$, we want to define $f^{\ast}([\alpha])$. The birational map $p^{-1} \circ f : \mathbb{P}^{2} \dashrightarrow X$ can be decomposed as a sequence of blow-ups and blow-downs. Hence we can choose a rational surface $(Y,q)$ in $\mathscr{E}$ such that the rational map $g:=p^{-1}\circ f\circ q$ is a morphism and define $f^{\ast}([\alpha]):=[g^{\ast}\alpha]$. 

\centerline{
\xymatrix{
{Y} \ar[r]^{g} \ar[d] & {X} \ar[d] \\
{\mathbb{P}^{2}} \ar[r]^{f} & {\mathbb{P}^{2}} \\
}
}

\noindent This defines a linear isometry $f^{\ast} : (\ZZ, (\cdot , \cdot)_{\mathbb{P}^{2}}) \to (\ZZ, (\cdot , \cdot)_{\mathbb{P}^{2}})$. We then  put $$f_{\ast}:=(f^{\ast})^{-1}$$ 
\noindent and observe that $f\mapsto f_{\ast}$ is a group homomorphism.

 We now describe a basis of the space $\ZZ$. On the union 
 $$\bigcup_{(X,p)\in \mathscr{E}}X$$
 
 \noindent of all rational surfaces in $\mathscr{E}$, we define an equivalence relation as follows: we say that a point $x$ in $(X,p)$ is equivalent to a point $y$ in $(Y,q)$ if the map $q^{-1}\circ p$ is well defined at $x$, maps $x$ onto $y$ and is a local isomorphism near $x$. We denote by $Ec(\mathbb{P}^{2})$ the quotient of $\cup_{(X,p)\in \mathscr{E}}X$ by this equivalence relation. If $x$ is a point in a rational surface $X$ dominating $\mathbb{P}^{2}$, we let $X_{x}\to X$ be the blow-up of $X$ at $x$ and $e_{x}\in H^{2}(X_{x},\Z)$ be the class of the exceptional divisor. We still denote by $e_{x}$ the image of this class in $\ZZ$ and observe that $e_{x}\in \ZZ$ depends only on the equivalence class of $x$ in $Ec(\mathbb{P}^{2})$. Hence for each point $a\in Ec(\mathbb{P}^{2})$ we have a class $e_{a}\in \ZZ$ associated to $a$. We have a natural map from the abelian group
  $$A(\mathbb{P}^{2}):=\Z\oplus \left( \underset{a\in Ec(\mathbb{P}^{2})}{\oplus}\Z_{a}\right)$$
  \noindent to $\ZZ$: it takes the first factor isomorphically to $H^{2}(\mathbb{P}^{2},\Z)\hookrightarrow \ZZ$ and the factor $\Z_{a}$ to $\Z e_{a}\subset \ZZ$. The direct sum above can be endowed with the bilinear form $(\cdot, \cdot)$ for which the vector $(t, (\lambda_{a})_{a\in Ec(\mathbb{P}^{2})})$ has norm $t^{2}-\sum_{a}\lambda_{a}^{2}$. We then have the following result: 
  
 \begin{prop} \cite{manin} The natural map
 $$(A(\mathbb{P}^{2}),(\cdot,\cdot)) \to (\ZZ,(\cdot ,\cdot)_{\mathbb{P}^{2}})$$
 \noindent is a surjective isometry. 
 \end{prop}

We will use two consequences of this proposition (see \cite{bfj,cantat}):
\begin{enumerate} 
\item First, one can complete $\ZZ$ to obtain a Hilbert space $\ZZZ$ endowed with a bilinear form (still denoted by $(\cdot , \cdot)_{\mathbb{P}^{2}}$) of signature $(1,\infty)$. Any isometry of $\ZZ$ extends to an isometry of $\ZZZ$, hence we get an action of the group $\bir$ on $\ZZZ$.  
\item Second, one can describe the action of a birational map $f$ on $\ZZ$ in terms of the basis given by the previous proposition. If $x\in \mathbb{P}^{2}$ and if $f$ is defined at $x$ we have: $$f_{\ast}(e_{x})=e_{f(x)}$$ \noindent (when one identifies $x$ and $f(x)$ to their images in $Ec(\mathbb{P}^{2})$). If $x$ is a point of a rational surface $X$ dominating $\mathbb{P}^{2}$, and if the image of $x$ in $\mathbb{P}^{2}$ is not a point of indeterminacy of $f$, the same formula remains true, provided that we define $f(x)$ correctly (see \cite{cantat}). In particular if $f$ is an automorphism of $\mathbb{P}^{2}$, $f_{\ast}$ permutes the vectors of the previous basis of $\ZZZ$. 
\end{enumerate}
We will denote by $u_{0}\in \ZZ$ the class of a line in $\mathbb{P}^{2}$. One can finally define the space $\mathbb{H}_{\mathbb{P}^{2}}$. This is simply the hyperbolic space associated to the space $(\ZZZ,(\cdot , \cdot)_{\mathbb{P}^{2}})$, i.e.:
 $$\mathbb{H}_{\mathbb{P}^{2}}:=\{u\in \ZZZ, (u,u)_{\mathbb{P}^{2}}=1, (u,u_{0})_{\mathbb{P}^{2}}>0\}.$$ 

\noindent {\bf Example 1, continued} {\it We denote by $x_{1}, x_{2}, x_{3}$ the three points $[1:0:0]$, $[0:1:0]$ and $[0:0:1]$ in $\mathbb{P}^{2}$. The group $A:=\C^{\ast}\times \C^{\ast}$ can be identified with the diagonal subgroup of ${\rm PGL}(3,\C)$. It fixes the class $u_{0}$ of a line as well as the three points $e_{x_{1}}$, $e_{x_{2}}$ and $e_{x_{3}}$ in $\ZZZ$. Hence it fixes pointwise the intersection of the $4$-dimensional subspace spanned by $u_{0}$ and the $e_{x_{i}}$'s with $\mathbb{H}_{\mathbb{P}^{2}}$, which is a copy of the $3$-dimensional hyperbolic space. 

To prove that the space $\left(\mathbb{H}_{\mathbb{P}^{2}}\right)^{A}$ is actually infinite dimensional, we observe that when we blow up the three points $x_{i}$ on $\mathbb{P}^{2}$, we obtain a rational surface $X$ on which $A$ acts with $6$ fixed points. By blowing up these six fixed points and iterating this procedure, we obtain a sequence of actions of $A$ on some rational surfaces $X_{n}$ dominating $\mathbb{P}^{2}$, for which the second Betti numbers $b_{2}(X_{n})$ go to infinity. This proves that the space of fixed points of $A$ in $\mathbb{H}_{\mathbb{P}^{2}}$ is infinite dimensional.}

We conclude this paragraph with a remark, which will be used below. If \mbox{$q : Y\to \mathbb{P}^{2}$} is a rational surface dominating $\mathbb{P}^{2}$, one can repeat everything we have said in this paragraph, replacing $\mathbb{P}^{2}$ by $Y$: one can define the spaces ${\rm Z}(Y)$, $\mathscr{Z}(Y)$, the form $(\cdot , \cdot)_{Y}$. One difference is that when defining the group $A(Y)$, the first $\Z$ factor has to be replaced by the group $H^{2}(Y,\Z)$. The spaces ${\rm Z}(Y)$ and $\mathscr{Z}(Y)$ are isomorphic to $\ZZ$ and $\ZZZ$, and when one identifies $\bir$ to ${\rm Bir}(Y)$ via $q$, the two actions we get on $\ZZZ$ and $\mathscr{Z}(Y)$ are conjugated. In particular one gets the following conclusion: if a birational map \mbox{$f : \mathbb{P}^{2}\dashrightarrow \mathbb{P}^{2}$} lifts to an automorphism of $Y$ which is isotopic to the identity, there is an orthonormal basis of $\ZZZ$ which is permuted by $f_{\ast}$. 

\subsection{Proof of Theorem~\ref{cremo}}
We now prove Theorem~\ref{cremo}. Let $\rho : \Gamma \to \bir$ be a non-elementary homomorphism. Let $\mathbb{H}_{\rho}$ be the unique minimal closed totally geodesic $\rho$-invariant subspace of $\mathbb{H}_{\mathbb{P}^{2}}$. We know that the action of $\Gamma$ on $\mathbb{H}_{\rho}$ factors through a homomorphism $\theta$ onto a subgroup $\Lambda$ of $\PSLR$. We need to show that $\Lambda$ is discrete.

We assume by contradiction that $\Lambda$ is dense in $\PSLR$. Let $\Lambda_{0}$ be a torsion-free finite index subgroup of $\Lambda$. The group $\Lambda_{0}$ is still dense in $\PSLR$. Since elliptic elements form an open set of $\PSLR$, there exists an element $\gamma_{0}\in \Lambda_{0}$ which is elliptic of infinite order. Let $g$ be an element of $\Gamma$ such that $\theta(g)=\gamma_{0}$. Then the birational transformation $\rho(g)$ has the following property: there exists a $2$-dimensional subspace $P\subset \ZZZ$ on which the form $(\cdot , \cdot)_{\mathbb{P}^{2}}$ is positive definite and on which $\rho(g)_{\ast}$ acts as an irrational rotation. Indeed, the harmonic map $f : \Delta \to \mathbb{H}_{\rho}\subset \mathbb{H}_{\mathbb{P}^{2}}$ is $\PSLR$-equivariant in our situation, hence has rank $2$ at every point of $\Delta$. One can thus take $P$ to be the image of the differential of $f$ at the point $\xi\in \Delta$ fixed by $\gamma_{0}$.

 The next proposition proves that such a behavior does not occur inside the Cremona group. Hence $\Lambda$ cannot be dense and has to be discrete. After eventually applying the Stein factorization theorem, we obtain a fibration of $X$ onto a hyperbolic $2$-orbifold $\Sigma$ and the kernel $H$ of the homomorphism $\Gamma \to \pi_{1}^{orb}(\Sigma)$ has to fix pointwise the subspace $\mathbb{H}_{\rho}$. Since $\mathbb{H}_{\rho}$ contains the image of $f$, its dimension is at least $2$. This concludes the proof of Theorem~\ref{cremo}.    

\begin{prop} Let $h : \mathbb{P}^{2} \dashrightarrow \mathbb{P}^{2}$ be a birational map such that the isometry $h_{\ast}$ fixes a point in $\mathbb{H}_{\mathbb{P}^{2}}$. Then there is no $h_{\ast}$-invariant $2$-dimensional subspace $P\subset \ZZZ$ on which $h_{\ast}$ acts as an irrational rotation.  
\end{prop}

{\it Proof.} According to \cite{cantat}, and up to replacing $h$ by a power, one can find a birational model $Y$ of $\mathbb{P}^{2}$ on which $h$ acts by an automorphism isotopic to the identity. The space $\mathscr{Z}(Y)$ is the direct sum of $H^{2}(Y,\R)$ and of $\ell^{2}(Ec(Y))$; $h_{\ast}$ acts trivially on the first factor and permutes the coordinates in the second. The plane $P$ has to be contained in the second factor. Let $u_{i}=\sum_{x\in Ec(Y)}\lambda_{i,x}e_{x}$ ($i=1,2$) be a basis of $P$. We know that there exists an angle $\alpha$ (with $\frac{\alpha}{\pi}$ irrational) such that:
$$\begin{array}{rcl}
h_{\ast}u_{1} & = & \cos(\alpha)u_{1}+\sin(\alpha)u_{2}, \\
h_{\ast}u_{2} & = & -\sin(\alpha)u_{1}+\cos (\alpha)u_{2}.\\
\end{array}$$
\noindent We fix a point $x\in Ec(Y)$ such that $(\lambda_{1,x},\lambda_{2,x})\neq (0,0)$. Looking coordinate by coordinate we obtain, for each integer $n$, the equalities:
$$\begin{array}{rcl}
\lambda_{1,h^{n}(x)} & = & \cos(n\alpha)\lambda_{1,x}+\sin(n\alpha)\lambda_{2,x},\\
\lambda_{2,h^{n}(x)} & = & -\sin(n\alpha)\lambda_{1,x}+\cos(n\alpha)\lambda_{2,x}.\\
\end{array}$$
Since the sum $\sum_{y}\lambda_{i,y}^{2}$ is finite, there must exist an integer $N$ such that $h^{N}(x)=x$. But this implies that the vector $(\lambda_{1,x},\lambda_{2,x})\in \R^{2}$ is an eigenvector for the rotation of angle $N\alpha$ in the plane. This is a contradiction.\hfill $\Box$

Theorem~\ref{cremo} only deals with non-elementary homomorphisms into the Cremona group but the work of Cantat also gives informations about elementary homomorphisms from finitely generated groups into $\bir$. In particular we have the following consequence of Theorem~\ref{serge-elliptic}.

\begin{center}
 {\it If $N$ is a finitely generated subgroup of the Cremona group which preserves a geodesic $\gamma$ in $\mathbb{H}_{\mathbb{P}^{2}}$ but is not elliptic, then either $N$ is virtually cyclic or $N$ is conjugated to a subgroup of the group $G_{toric}$.}
 \end{center}
 
\noindent Let us explain this fact. Let $u : N \to \R$ be the ``translation homomorphism" induced by the action of $N$ on $\gamma$. It is known that its image is cyclic, see \cite{cantat,favre}; this uses the fact that $N$ is finitely generated. If $N$ is not virtually cyclic, the group $A:={\rm Ker}\, u$ is infinite and elliptic and we are in the situation of Theorem~\ref{serge-elliptic}, establishing the above result.

We can now prove Corollary~\ref{sc}. 

{\it Proof of Corollary~\ref{sc}.} Let $\rho : \Gamma_{1} \to \bir$ be an injective homomorphism where $\Gamma_{1}$ is a cocompact lattice in $\SU$ ($n\ge 2$). Note that if the conclusion of the corollary holds for a normal subgroup of finite index of $\Gamma_{1}$, it holds also for $\Gamma_{1}$. Hence we can assume that $\Gamma_{1}$ is torsion-free. We assume by contradiction that $\rho$ is non-elementary, and apply Corollary~\ref{nonelecomplet}: $\Gamma_{1}$ has no finite index subgroup isomorphic to a surface group, hence $\Gamma_{1}$ must be isomorphic to a subgroup of the group $G_{toric}$. But $\Gamma_{1}$ has no nontrivial normal abelian subgroup, hence its intersection with $\C^{\ast}\times \C^{\ast}\subset G_{toric}$ must be trivial. This implies that $\Gamma_{1}$ is isomorphic to a subgroup of ${\rm GL}(2,\Z)$, a contradiction. Hence the homomorphism $\rho$ has to be elementary. If $\rho(\Gamma_{1})$ preserves a geodesic we apply the remark made above. Since $\Gamma_{1}$ is not virtually cyclic, it has to be isomorphic to a subgroup of the group $G_{toric}$, yielding a contradiction again.  

Finally, if $\rho$ is elementary but does not preserve a geodesic, it has to fix a point in the space $\mathbb{H}_{\mathbb{P}^{2}}$ or to fix a unique point of the boundary of $\mathbb{H}_{\mathbb{P}^{2}}$.\hfill $\Box$

%%%%%%%%%%%%%%%%%%%%%%%%%%%%%%%%%%%%%%%%%%%%%%%%%%%%%%%%%%%%%%%%%%%%%%%%%%%%%%%%%%%%%%%
%%%%%%%%%%%%%%%%%%%%%%%%%%%%%%%%%%%%%%%%%%%%%%%%%%%%%%%%%%%%%%%%%%%%%%%%%%%%%%%%%%%%%%%

\section{Koszul and Malgrange's integrability theorem}\label{section-dbar}

In this section, we prove the infinite dimensional version of Koszul and Malgrange's integrability theorem.

Let us first remind the content of the proof in the finite dimensional case. Assume that $E\to X$ is a smooth complex vector bundle over the complex manifold $X$ and that $D$ is a $(0,1)$-connexion on $E$, i.e. a differential operator $D : A^{0}(E)\to A^{1}(E)$ (where $A^{0}(E)$ is the space of smooth sections of $E$ and $A^{1}(E)$ the space of smooth $1$-forms on $X$ with values in $E$) such that:
$$D(fs)=\overline{\partial}(f )s+fDs$$

\noindent whenever $f$ is a smooth function and $s$ is a section of $E$. The operator $D$ naturally extends to an operator acting on all differential forms with values in $E$. The theorem of Koszul and Malgrange asserts that if $D$ satisfies the integrability condition $D^{2}=0$, one can find a holomorphic structure on $E$ for which $D$ is the usual $\overline{\partial}$ operator. This is a simple particular case of the integrability criterion of Newlander and Nirenberg. The proof can be described as follows. Let $s_{1},\ldots ,s_{n}$ be any local frame for the bundle $E\to X$, over some open set $U\subset X$. We want to prove that around any point, one can find a frame made of sections annihilated by $D$ (the transition matrices between any to such frames will then be holomorphic). We can write:

$$Ds_{j}=\sum_{i=1}^{n}a_{ij}s_{i},$$

\noindent where $A=(a_{ij})$ is a matrix of $(0,1)$-forms on $U$. We look for an invertible matrix-valued map $M=(m_{ij})$ on $U$ for which the sections $u_{j}=\sum_{i=1}^{n}m_{ij}s_{i}$ would satisfy $Du_{j}=0$.
In other words we want to solve (locally) the equation $\overline{\partial} M=-AM$, knowing the integrability condition
$$\overline{\partial}A+A\wedge A=0,$$

\noindent which is equivalent to the condition $D^{2}=0$. The (matrix valued) $2$-form  $A\wedge A$ is defined by:
$$A\wedge A =\sum_{j<l}[A_{j},A_{l}]d\overline{z}_{j}\wedge d\overline{z}_{l},$$
\noindent where $A=\sum_{i=1}^{k}A_{i}d\overline{z}_{i}$. Here $(z_{1}, \ldots ,z_{k})$ are coordinates on $U$ and the $A_{i}=A_{i}(z)$ are $n\times n$ matrices.

 In the infinite dimensional setting, $A=\sum_{i=1}^{k}A_{i}d\overline{z}_{i}$ where the $A_{i}$'s are smooth maps with values in ${\rm End}(\mathscr{H})$ ($\mathscr{H}$ being a fixed complex Hilbert space) and we look for a smooth map $M$ with values in the open set $\GL\subset {\rm End}(\mathscr{H})$ of invertible endomorphisms of $\mathscr{H}$, solving the equation above. To solve this problem, we use the following:

\begin{prop}\label{edp} Let $C$ be a smooth map with values in the space ${\rm End}(\mathscr{H})$, defined in a neighborhood of $0$ in $\C^{k}$. Assume that $C$ is holomorphic with respect to $z_{1},\ldots , z_{p}$ for some integer $0\le p < k$. Then there exists a smooth map $h$ defined in a neighborhhood of $0$ in $\C^{k}$, with values in the space $\GL$ and holomorphic in $z_{1},\ldots , z_{p}$, such that:
$$h^{-1}\frac{\partial h}{\partial \overline{z}_{p+1}}=C.$$\end{prop}

To deduce the theorem from the proposition, one proceeds exactly as in \cite{km}. Note that it is at this point that we use the integrability condition satisfied by $A$. Since this part of the proof is exactly the same as in the finite dimensional case, we refer the reader to \cite{km}, page 103. We now turn to the proof of Proposition~\ref{edp}. 

We follow the proof explained in~\cite{donald}, page 50. We denote by $B$ a ball in $\C$, by $E$ a complex Banach space, and by $C^{l+\alpha}(B,E)$ the usual H\"{o}lder space of $C^{l}$ maps from $B$ to $E$, with $\alpha$-H\"{o}lder derivatives of order $l$ ($0< \alpha <1$). We will use the following result:

\begin{prop} There exists a continuous linear operator $P : C^{\alpha}(B,E)\to C^{1+\alpha}(B,E)$ such that: $$\frac{\partial}{\partial \overline{z}}\left(Pg \right)=g,$$ 

\noindent for all $g\in C^{\alpha}(B,E)$. Moreover, the operator $P$ is continuous from $C^{l+\alpha}(B,E)$ to $C^{l+1+\alpha}(B,E)$ for every integer $l$. 
\end{prop}
{\it Proof.}  The construction of the operator $P$ is classical in the case of complex-valued functions, see~\cite{vekua} (chapter I, \S \S 6, 8, 9) for its construction and continuity between various H\"{o}lder and Sobolev spaces. Here we are dealing with functions with values in a complex Banach space. Observe that all the integrals of functions with values in a Banach space that we need to consider to define $P$ are integrals of {\it continuous functions}. They can thus be defined by limits of Riemann sums, and the separability (or not) of the space $E$ plays no role here. We shall indeed apply this proposition to some non separable Banach space, namely, the space of endomorphisms of a Hilbert space.

 Therefore one sees readily that the proofs given in \cite{vekua} extend to this context. We only give the formula defining $P$ and refer to \cite{vekua} for its properties:
$$P(g)(z)=\frac{1}{2i\pi}\int_{B}\frac{g(\xi)}{\xi -z}d\xi \wedge d\overline{\xi}.$$\hfill $\Box$

We now prove Proposition~\ref{edp}. Write $z=z_{p+1}$ and $z'=(z_{1},\ldots , z_{p},z_{p+2},\ldots, z_{k})$. We first treat $z'$ as a parameter. We thus want to solve the equation
$$\frac{\partial h}{\partial \overline{z}}=hC,$$

\noindent in a neighborhood of $0\in \C$. A simple renormalization argument (see~\cite{donald} page 51 for instance) implies that it is enough to solve this equation when $C$ is defined on the ball $B$ of radius $1$ around $0$ in $\C$ and the $C^{1}$ norm of $C$ is as small as we wish. By looking for $h$ under the form $h={\rm Id}+u$ and by using the operator $P$ given by the previous proposition (with the Banach space $E$ being the space of endomorphisms of $\mathscr{H}$), we see that it is enough to find a solution of the equation:
$$u=P(C+uC),$$
\noindent where $u\in C^{\alpha}(B,E)$. If $C$ is small enough we can apply the contraction mapping principle to the operator $u\mapsto P(C+uC)$ in the space $C^{\alpha}(B,E)$. Its fixed point will be close to $0$, hence the map $h={\rm Id}+u$ will take its values in ${\rm GL}(\mathscr{H})$.

The fact that the solution we obtain depends smoothly on $z'$ and holomorphically on the $z_{i}$ ($i\le p$) is a consequence of the implicit function theorem. This proves Proposition~\ref{edp}.

%%%%%%%%%%%%%%%%%%%%%%%%%%%%%%%%%%%%%%%%%%%%%%%%%%%%%%%%%%%%%%%%%%%%%%%%%%%%%%%%%%%%%%%%%%%%%%%%%%%%%%%%%%%%%%%%%%%%%%%%%%%%%%%%%%%%%%%%%%%%%%%%%%%%%%%%%%%%%%%%%%%%%%%%%%%%%%%%%%%%%%%%%%%%%%%%%%%%%%%%%%%%%%%%

\section{About the regularity of harmonic maps}\label{analyse}

The purpose of this section is to briefly explain the proof of Proposition~\ref{laplacelp}. We first state a version of Proposition~\ref{laplacelp} for the flat Laplacian $\Delta$, i.e. the Laplacian associated to the euclidean metric on $\C^{n}$. Here again $B$ is a closed ball in $\widetilde{X}$ that we can identify with a ball in $\C^{n}$ via some choice of coordinates. 

\begin{prop}\label{plat} Let $p>1$. There exists a constant $C$ such that if $u\in W^{1,p}(B,\mathscr{H})$ and $\Delta u\in L^{p}(B,\mathscr{H})$ then $u\in W^{2,p}(B,\mathscr{H})$ and we have: $\vert D^{2} u\vert_{p} \le C \vert \Delta u\vert_{p}$. Here $D^{2}u$ is the matrix of partial derivatives of order $2$ of $u$ and $\vert D^{2}u\vert_{p}^{p}=\sum_{i,j}\vert \partial_{i}\partial_{j}u\vert_{p}^{p}$.
\end{prop}

\noindent {\it Proof.} This is classical. It follows from the fact that for each index $\alpha$, there exists an linear operator $R_{\alpha}$ (the Riesz transform) mapping continuously $L^{p}(B,\mathscr{H})$ into itself for each $p\in (1,+\infty)$ and such that $\frac{\partial^{2}u}{\partial x_{\alpha}\partial x_{\beta}}=-R_{\alpha}R_{\beta}(\Delta u)$. See \cite{stein}, chapter 2 for the discussion of the continuity properties of various operators defined by singular integrals. The case of vector-valued functions is discussed in section 5 of Chapter 2. The particular case of the Riesz transforms is discussed in Chapter 3.\hfill $\Box$

With Proposition~\ref{plat} at hand, one proves Proposition~\ref{laplacelp} exactly as in the finite dimensional case. One first proves the following {\it a priori} estimate: 

\begin{prop}\label{nonplat} Let $B'\subset \subset B''$ be closed balls contained in the interior of $B$. There exists a constant $C$ such that if $u\in W^{2,p}_{loc}(B)$ then $$\vert D^{2}u\vert_{p,B'}\le C \vert \Delta_{g} u\vert_{p,B''}.$$
\end{prop}

\noindent From this statement, one easily deduces the result of Proposition~\ref{laplacelp}. We refer the reader to~\cite{gt}, Chapter 9 for more details.

\addcontentsline{toc}{section}{Appendix.-- \; Elliptic subgroups of the Cremona group with a large normalizer, by Serge Cantat}
\section*{Appendix.-- \; Elliptic subgroups of the Cremona group with a large normalizer, by Serge Cantat}\label{app}

In this appendix, we prove the following theorem (Theorem~\ref{serge-elliptic} in the introduction).  Its statement uses
the action of $\bir$ on $\mathbb{H}_{\mathbb{P}^{2}}$. Recall that the distance $d(u,v)$ on 
$\mathbb{H}_{\mathbb{P}^{2}}$ is given by the intersection form $(u,v)_{\mathbb{P}^{2}}$: $\cosh(d(u,v))=(u,v)_{\mathbb{P}^{2}}$. Birational transformations determine three kind of isometries of this space: elliptic, parabolic, or hyperbolic.

\setcounter{main}{0}
\setcounter{rem}{0}
\setcounter{prop}{0}

\begin{main}\label{thm:Serge}
Let $N$ be a subgroup of the Cremona group $\bir$. Assume that there exists a short exact
sequence 
\[
1 \to A \to N \to B\to 1
\]
where $N$ contains at least one hyperbolic element,  and $A$ is infinite elliptic. Then $N$ is conjugate to a subgroup of the group $G_{toric}$ of automorphisms of $(\C^*)^2$.
\end{main}

\begin{rem}
The same conclusion holds if we assume that  $A$ is infinite, finitely generated, and all its
elements are elliptic, because this new assumption implies that $A$ is elliptic (see Theorem 6.4 and Proposition 6.12 in \cite{cantat}).
\end{rem}

Let $Y$ be a rational surface. The group ${\rm Aut}(Y)$ of automorphisms of $Y$ is a complex Lie group. 
It may have an infinite number of connected components, but its connected component of the identity ${\rm Aut}(Y)^0$ is 
easily described: ${\rm Aut}(Y)^0$ is isomorphic to an algebraic subgroup of ${\rm PGL}_{n+1}(\C)$ for some positive integer $n$; this isomorphism is given by an equivariant embedding of $Y$ in $\mathbb{P}^{n}(\C)$. All pairs $(Y, {\rm Aut}(Y)^0)$ with $\dim({\rm Aut}(Y)^0)\geq 1$  are described in \cite{akhiezer, blanc}. 

Let $u_0\in \ZZZ$ be the class of a line in~$\mathbb{P}^2$. Since $A$ is elliptic, it fixes a point $u$ in $\mathbb{H}_{\mathbb{P}^{2}}$ and we have: $$(a_*(u_0), u_0)_{\mathbb{P}^{2}}\leq \cosh(2d(u_0,u)),$$ for all $a$ in $A$. Since $(a_*(u_0), u_0)_{\mathbb{P}^{2}}$ is the degree of $a$, $A$ is contained in the algebraic variety ${\rm Bir}_D(\mathbb{P}^2)$ of birational transformations of the plane of degree $\leq D$, for $D=[\cosh(2{\rm dist}(u_0,u))]+1$. Let $\overline{A}$ be the Zariski closure of $A$ in ${\rm Bir}_D(\mathbb{P}^2)$; since $A$ is infinite, $\dim({\overline{A}})\geq 1$; since $A$ is normalized by $N$, so is $\overline{A}$. From Weil and Rosenlicht Theorem (see \cite{rosenlicht,weil}, and \cite{blanc}), there exist a  surface $Y$ and a birational map $\pi:Y\dasharrow \mathbb{P}^2(\C)$ such that $\pi^{-1}\circ {\overline{A}}\circ \pi$ is contained in the group ${\rm Aut}(Y)$ and intersects  a finite number of its components.
After conjugation by $\pi$, we can and do assume that $N$ is contained in ${\rm Bir}(Y)$ and  $\overline{A}$ is a Zariski closed complex Lie subgroup of ${\rm Aut}(Y)$. Changing $A$ into $\overline{A}$ and changing $N$ accordingly, we have now a short exact sequence
\[
1 \to {\overline{A}} \to  {\overline{N}} \to {\overline{B}} \to 1
\]
in which $ {\overline{N}}$ contains hyperbolic elements and $ {\overline{A}}$ is a Lie subgroup of ${\rm Aut}(Y)$ with 
finitely many components. Let $ {\overline{A}}_0$ be the connected component of the identity in $ {\overline{A}}$. 

\begin{lemma}
The group $ {\overline{N}}$ normalizes $ {\overline{A}}_0$, and the group ${\overline{A}}_0$ has a Zariski open orbit in $Y$.
\end{lemma}

\noindent {\it Proof.} Since ${\overline{N}}$ normalizes ${\overline{A}}$, it normalizes the connected component $ {\overline{A}}_0$.
Since $ {\overline{A}}_0$ is Zariski-closed in the algebraic group ${\rm Aut}(Y)^0$, its orbits in $Y$ 
are Zariski open subsets of Zariski closed subsets in $Y$.  Thus, two cases can occur. Either
$ {\overline{A}}_0$ has an open orbit, or the orbits of $ {\overline{A}}_0$ determine a pencil 
of curves on $Y$. In the second case, the group $\overline{N}$ permutes the members of the 
pencil; in particular, all elements of $N$  preserve a meromorphic  fibration on $Y$, and $N$
does not contain any hyperbolic element (see \cite{cantatfavre}). From this contradiction, we
deduce that $ {\overline{A}}_0$ has an open orbit in $Y$; in particular, $\dim({\overline{A}}_0)\geq 2$.\hfill $\Box$

Let now $f$ be an element of $\overline{N}$. Since $f$ normalizes $ {\overline{A}}_0$ and $ {\overline{A}}_0$ 
is contained in ${\rm Aut}(Y)$, the group $ {\overline{A}}_0$  permutes the indeterminacy set ${\rm Ind}(f)$ and the
set of exceptional curves ${\rm Exc}(f)$ (${\rm Exc}(f)$ is the zero locus of the jacobian determinant of $f$). In particular, 
both ${\rm Ind}(f)$ and ${\rm Exc}(f)$ are contained in the complement of the open orbit of  $ {\overline{A}}_0$; this proves the following lemma.

\begin{lemma}
The union of all indeterminacy sets ${\rm Ind}(f)$  and exceptional sets ${\rm Exc}(f)$ of elements of $ {\overline{N}}$
is a proper Zariski closed set $Z\subset Y$. The group $\overline{N}$ acts by automorphisms on the
Zariski open set $U=Y\setminus Z$. 
\end{lemma}

We can now restrict our study to the Zariski open set $U\subset Y$. 
Each element $f$ of $\overline{N}$ determines an automorphism $\varphi_f$ of the algebraic
group ${\overline{A}}_0$ such that 
\[
f\circ a \circ f^{-1}=\varphi_f(a)
\] 
for all $a$ in ${\overline{A}}_0$. The group ${\overline{A}}_0$ splits into a $\varphi_f$-invariant exact sequence 
\[
1\to R \to {\overline{A}}_0 \to S \to 1
\]
where $R$ is the solvable radical of ${\overline{A}}_0$ and $S$ is semi-simple. 

Assume, first, that $R$ is trivial. Then ${\overline{A}}_0=S$ is semi-simple and the group of exterior automorphisms
of ${\overline{A}}_0$ is finite. Let $f$ be a hyperbolic element of ${\overline{N}}$. Changing $f$ into one of its positive iterates, 
$\varphi_f$ is an interior automorphism: there exists $b$ in $S$ such that $f\circ a \circ f^{-1}=b\circ a \circ b^{-1}$ for all $a$ in $S$. 
In particular, $f$ commutes to $b$; since $b$ is in ${\rm Aut}^0(Y)$ and $f$ is hyperbolic, we deduce from \cite{cantat}, Theorem B, that
$b$ has finite order. Changing $f$ into a new positive iterate, we can assume that $b$ is the identity and thus $f$ commutes
to the action of $S$. Theorem B of \cite{cantat} provides  a contradiction because $f$ is hyperbolic. This shows that {\sl{the solvable radical of ${\overline{A}}_0$ is non trivial}}.

The radical $R$, and all members
of its central series $R_i=[R_{i-1},R_{i-1}]$ are 
$\varphi_f$-invariant. Thus, changing ${\overline{A}}_0$ into the last non trivial derived subgroup ${\overline{A}}_1$ of $R$, we 
have a new exact sequence 
$1 \to {\overline{A}}_1 \to {\overline{N}} \to Q \to 1
$ where ${\overline{A}}_1$ is a connected abelian group with an open orbit $V\subset U$. 

There are three possibilities for this algebraic group ${\overline{A}}_1$: $\C^2$, $\C\times \C^*$, and $\C^*\times \C^*$. 
To conclude, we need to exclude the first two cases. 
In the first case, ${\overline{A}}_1= \C^2$, the unique open orbit
of ${\overline{A}}_1$ in $U$ can be identified to $\C^2$ and all elements $f$ of ${\overline{N}}$ act by affine automorphisms
on $V=\C^2$ (with linear part given by $\varphi_f\in {\rm GL}(2,\C)$). In particular, the degree  of $f^n$ is bounded and
${\overline{G}}$ does not contain any hyperbolic element, contradicting our assumption. 
In the second case, ${\overline{A}}_1= \C \times \C^*$. The group of algebraic automorphisms of ${\overline{A}}_1$ 
is made of diagonal transformations $(x,y)\mapsto (ax,\epsilon y)$ with $a \in \C^*$ and  $\epsilon =\pm 1$; thus, the
same argument applies to exclude this case.   

As a consequence, ${\overline{A}}_1=\C^*\times \C^*$ and there is a ${\overline{N}}$-invariant Zariski-open subset $V\subset Y$ such that  $V$ is isomorphic to ${\overline{A}}_1$ and ${\overline{N}}$ acts by regular transformations on $V$. Since the group of regular
automorphisms of $\C^*\times \C^*$ is the semi-direct product of $\C^*\times \C^*$ (acting by translations) and ${\rm GL}(2,\Z)$
(acting by group automorphisms, i.e. by monomial transformations), Theorem \ref{thm:Serge} is proved.

\bigskip
\bigskip
\begin{small}
\begin{tabular}{llll}
Thomas Delzant & & & Pierre Py\\
IRMA & & & Department of Mathematics\\
Universit\'e de Strasbourg \& CNRS & & & University of Chicago\\
67084 Strasbourg, France & & & Chicago, Il 60637, USA\\
delzant@math.unistra.fr & & & pierre.py@math.uchicago.edu\\    
\end{tabular}

\medskip

\begin{tabular}{ll}
Serge Cantat & \\
IRMAR (UMR 6625 du CNRS) & \\
Universit\'e de Rennes 1 & \\
35042 Rennes, France & \\
 serge.cantat@univ-rennes1.fr & \\
\end{tabular}
\end{small}

 \end{document}